\documentclass[3p,sort&compress]{elsarticle}
\journal{Journal of Computational and Applied Mathematics (Elsevier)}

\usepackage{framed,multirow}
\usepackage{latexsym}
\usepackage{amsmath}
\usepackage{amssymb}
\usepackage{amsthm}
\usepackage{mathtools}
\usepackage{pifont}
\usepackage{enumitem}
\usepackage{booktabs}
\usepackage{subcaption}
\usepackage[dvipsnames]{xcolor}
\usepackage{yfonts}
\usepackage[scr=boondoxo,scrscaled=1.05]{mathalfa}
\usepackage{appendix}
\usepackage{enumitem}
\usepackage{url}
\usepackage{ifthen}
\usepackage{tensor}

\theoremstyle{plain}

\theoremstyle{definition}
\newtheorem{definition}{Definition}

\theoremstyle{remark}
\newtheorem{remark}{Remark}

\begin{document}

\begin{frontmatter}
\title{Flow-driven spectral chaos (FSC) method for long-time integration of second-order stochastic dynamical systems}

\author[1]{Hugo Esquivel}
\ead{hesquive@purdue.edu}

\author[1]{Arun Prakash}
\ead{aprakas@purdue.edu}

\author[2]{Guang Lin\fnref{fn2}}
\ead{guanglin@purdue.edu}

\address[1]{Lyles School of Civil Engineering, Purdue University, 550 W Stadium Ave, West Lafayette IN 47907, USA}
\address[2]{Department of Mathematics, School of Mechanical Engineering, Purdue University, 150 N University St, West Lafayette IN 47907, USA}

\fntext[fn2]{Department of Statistics, Department of Earth, Atmospheric, and Planetary Sciences (by courtesy), Purdue University.}

\begin{abstract}
For decades, uncertainty quantification techniques based on the spectral approach have been demonstrated to be computationally more efficient than the Monte Carlo method for a wide variety of problems, particularly when the dimensionality of the probability space is relatively low.
The time-dependent generalized polynomial chaos (TD-gPC) is one such technique that uses an evolving orthogonal basis to better represent the stochastic part of the solution space in time.
In this paper, we present a new numerical method that uses the concept of \emph{enriched stochastic flow maps} to track the evolution of the stochastic part of the solution space in time.
The computational cost of this proposed flow-driven stochastic chaos (FSC) method is an order of magnitude lower than TD-gPC for comparable solution accuracy.
This gain in computational cost is realized because, unlike most existing methods, the number of basis vectors required to track the stochastic part of the solution space, and consequently the computational cost associated with the solution of the resulting system of equations, does not depend upon the dimensionality of the probability space.
Four representative numerical examples are presented to demonstrate the performance of the FSC method for long-time integration of second-order stochastic dynamical systems in the context of stochastic dynamics of structures.
\end{abstract}

\end{frontmatter}

\section*{Highlights}

\begin{itemize}\setlength\itemsep{0em}
\item A novel and efficient numerical method that uses an enriched stochastic flow map to track the evolution of the solution space is proposed for long-time integration of stochastic dynamics of structures.
\item For the same level of accuracy, the computational cost associated with the proposed FSC method is an order of magnitude lower than other time-dependent, spectral-based methods in the current literature.
\item It is partially insensitive to the curse of dimensionality compared to other spectral-based methods such as gPC and TD-gPC.
\item It has the potential to be used in the context of large-scale structural engineering problems to quantify uncertainties of long-time response with high fidelity.
\item Exact response expressions for the mean and variance of a single-degree-of-freedom system under free vibration and uniform stochastic stiffness are also provided.
\end{itemize}

\section{Introduction}\label{sec1Introduction}

In the past few decades, the area of uncertainty quantification has received increasing attention, primarily in the fields of physics and engineering.
This is not surprising since any mathematical description of a real-life physical system is always subject to the effects of input randomness.
In structural engineering, for example, the randomness of a system may arise from different sources, such as variability of material properties, imperfections in geometry, loading scenarios, boundary conditions, etc.
Though these physical quantities may be random, in most cases their stochastic characteristics can be identified and modeled mathematically using one or more of the numerous distributions available in statistics.

Due to the inherent complexity of real-life systems, closed-form solutions are not always possible.
As an alternative, one uses accurate numerical methods that allow one to propagate and quantify the effects of input uncertainties in the system response efficiently.
Historically, the polynomial chaos (PC) method developed by Wiener \cite{wiener1938homogeneous} and extended by Ghanem and Spanos \cite{ghanem1990polynomial} in the context of stochastic finite elements, has been used for uncertainty quantification.
Xiu and Karniadakis \cite{xiu2002wiener} further developed the generalized-PC (gPC) method in which time-independent polynomials are used to decompose a stochastic process into deterministic and non-deterministic parts using the orthogonality property of the basis in the random function space (RFS).
The benefit of using such orthogonal polynomials in the RFS is that when the underlying process (i.e.~the solution of the ODE) is represented with them, the method leads to exponential convergence to the solution (provided the stochastic part of the solution space is not discretized).

The gPC method has undergone several modifications and extensions to improve its computational efficiency, its effectiveness for long-time integration, and its ability to handle stochastic discontinuities.
The \emph{multi-element gPC} and related methods, developed by Karniadakis and others \cite{wan2005adaptive,wan2006multi,wan2006long,foo2008multi,agarwal2009domain,ma2009adaptive,kewlani2009multi,foo2010multi,prempraneerach2010uncertainty,jakeman2013minimal}, adaptively decompose the probability space into elements until a pre-specified threshold for the relative error in variance is reached.
Then a stochastic spectral expansion on each element is used to advance the state of the system one-time step forward, and this process is repeated every time the threshold is exceeded during the simulation.
Another approach is the \emph{dynamically orthogonal PC} (DO-PC) method \cite{sapsis2009dynamically,choi2013convergence,ueckermann2013numerical,cheng2013dynamicallyI,cheng2013dynamicallyII} where the time rate of change of the spatio-temporal function space is ensured to be orthogonal to itself.
This condition, called the dynamically orthogonal (DO) condition, is enforced at every time step as the simulation proceeds.
The DO-PC scheme is essentially a generalization of the POD (Proper Orthogonal Decomposition) method \cite{papoulis1965random,holmes2012turbulence} and the gPC method in the framework of continuous stochastic dynamical systems.
An error analysis for the DO-PC method can be found in \cite{musharbash2015error}.

For long-time dynamical simulations, when the stochasticity of the system has developed significantly, the gPC method fails to capture the probability moments accurately because the probability distribution of the solution changes significantly with time.
The time-dependent gPC (TD-gPC) method was proposed by Gerritsma et al.~\cite{gerritsma2010time} to allow the basis to evolve in time so as to better represent the transient nature of the probability distribution of the solution during the simulation.
Heuveline and Schick \cite{heuveline2014hybrid} modified the TD-gPC method (mTD-gPC) to account for systems governed by second-order ODEs and also improved the accuracy of the method.
Generally speaking, the mTD-gPC method relies on spanning the stochastic part of the solution space, at certain time steps (aka \emph{reset times}), by taking a full tensor product between an evolving RFS (that depends upon the evolution of the state variables of the system) and the original RFS.
However, since a full tensor product is required to be conducted at the reset times, the method suffers from the curse of dimensionality.
A \emph{hybrid generalized polynomial chaos} was also developed in \cite{heuveline2014hybrid} to help alleviate the curse of dimensionality of mTD-gPC.
This method, however, requires the use of an $(h,p)$-discretization for the stochastic part of the solution space in contrast to mTD-gPC which solely requires the use of a $p$-discretization.

More recently, a generalized PC method based on flow map composition was proposed by Luchtenburg et al.~\cite{luchtenburg2014long} to address the long-time uncertainty propagation in dynamical systems more effectively.
The method fundamentally uses short-time flow maps based on spectral polynomial bases to account for the stretching and folding effect caused by the evolution of the system's state over time.
As with mTD-gPC, this method suffers from the curse of dimensionality, because a tensor product is needed to construct the spectral polynomial bases during the simulation.
Ozen and Bal \cite{ozen2016dynamical} developed the \emph{dynamical gPC} (DgPC) method to address the long-time behavior of stochastic dynamical systems via a generalization of the PCE (Polynomial Chaos Expansion) framework. 
They demonstrated that results from DgPC match well with those obtained from other standard methods such as Monte Carlo.

In this paper, a flow-driven spectral chaos (FSC) method is proposed which tracks the evolution of the random basis via an \emph{enriched stochastic flow map} of the state of the system.\footnote{A \emph{generalized} version of this FSC method can already be found in \cite{esquivel2020flow}. In that work, the FSC scheme presented in Section \ref{sec1FSCSch} is simply denoted as FSC-1.}
The enriched flow map consists of the time derivatives of the solution up to a specified order.
Unlike any gPC-based method, the number of basis vectors needed in FSC to represent the stochastic part of the solution space does not grow with the dimensionality of the probability space.
However, as with any gPC-based method, the FSC method requires the computation of inner products using quadrature points distributed over the entire random domain, and the number of quadrature points can grow significantly with the dimensionality of the probability space.
Nevertheless, FSC presents a significant advantage over gPC-based methods as its computational cost, for comparable accuracy, is an order of magnitude lower compared to existing TD-gPC methods.
Conversely, for comparable computational cost, FSC is able to achieve an order of magnitude of better accuracy than TD-gPC.

This paper is organized as follows.
Section \ref{sec1SettingNotation} introduces the notation and definitions used in this paper and Section \ref{sec1ProSta} presents the precise problem statement.
The spectral approach for solving this problem is presented in Section \ref{sec1SolSpeApp} and Section \ref{sec1FSCMet} describes the proposed FSC method in detail.
Three numerical examples are presented in Section \ref{sec1NumRes} to demonstrate and compare the accuracy of FSC to other existing methods, such as mTD-gPC and Monte Carlo.
The FSC method is then applied to quantify uncertainties in the structural dynamics of a 3-story building subject to an earthquake excitation in Section \ref{sec1AppStrDyn}.
In Appendix \ref{appsec1RanBasTimInt}, we define the random bases that we use to span the random function space of the problem in hand.
Finally, in Appendix \ref{appsec1SDOF}, we provide expressions for the mean and variance of the exact response of a single-degree-of-freedom system under free vibration and uniformly-distributed stochastic stiffness, followed by a brief outline of the standard Monte Carlo method in Appendix \ref{appsec1OveStaMonCarMet}.

\section{Setting and notation}\label{sec1SettingNotation}

\begin{definition}[Temporal space]
Let $(\mathfrak{T},\mathcal{O})$ be a topological space, where $\mathfrak{T}=[0,T]$ is a closed interval representing the temporal domain of the system, $T$ is a positive real number symbolizing the duration of the simulation, and $\mathcal{O}=\mathcal{O}_\mathbb{R}\cap\mathfrak{T}$ is the topology on $\mathfrak{T}$ with $\mathcal{O}_\mathbb{R}$ denoting the standard topology over $\mathbb{R}$.
In this paper, $(\mathfrak{T},\mathcal{O})$ is called the \emph{temporal space} of the system.
\end{definition}

\begin{definition}[Random space]
Let $(\Omega,\boldsymbol{\Omega},\lambda)$ be a (complete) probability space, where $\Omega$ is the sample space, $\boldsymbol{\Omega}\subset 2^\Omega$ is the $\sigma$-algebra on $\Omega$ (aka the collection of events in probability theory), and $\lambda:\boldsymbol{\Omega}\to[0,1]$ is the probability measure on $\boldsymbol{\Omega}$.
Let $\xi:(\Omega,\boldsymbol{\Omega})\to(\mathbb{R}^d,\mathcal{B}_{\mathbb{R}^d})$ be a measurable function (aka random variable) given by $\xi=\xi(\omega)$, with $\mathcal{B}_{\mathbb{R}^d}$ denoting the Borel $\sigma$-algebra over $\mathbb{R}^d$.
In this work, the measure space $(\Xi,\boldsymbol{\Xi},\mu)$ is called \emph{random space}, where $\Xi=\xi(\Omega)\subset\mathbb{R}^d$ is a set representing the random domain of the system, $d$ denotes the dimensionality of the random space, $\boldsymbol{\Xi}=\mathcal{B}_{\mathbb{R}^d}\cap\Xi$ is the $\sigma$-algebra on $\Xi$, and $\mu:\boldsymbol{\Xi}\to[0,1]$ is the probability measure on $\boldsymbol{\Xi}$ defined by the pushforward of $\lambda$ by $\xi$, that is $\mu=\xi_*(\lambda)$.
\end{definition}

Note that more structure can be added to these spaces whenever they are needed in the analysis.
However, in order to keep the above definitions as elementary as possible, we singled out those mathematical objects that did not play a crucial role in the development of this work, such as the definition of a metric, a norm or an inner product for the underlying set.

\begin{definition}[Temporal function space]
Let $\mathscr{T}(n)=C^n(\mathfrak{T},\mathcal{O};\mathbb{R})$ be a continuous $n$-differentiable function space.
This \emph{temporal function space} is the space of all functions $f:(\mathfrak{T},\mathcal{O})\to(\mathbb{R},\mathcal{O}_\mathbb{R})$ that have continuous first $n$ derivatives on $(\mathfrak{T},\mathcal{O})$.
\end{definition}

\begin{definition}[Random function space]\label{sec1defRanFunSpa}
Let $\mathscr{Z}=(L^2(\Xi,\boldsymbol{\Xi},\mu;\mathbb{R}),\langle\,\cdot\,,\cdot\,\rangle)$ be a Lebesgue square-integrable space equipped with its standard inner product
\begin{equation*}
\langle\,\cdot\,,\cdot\,\rangle:L^2(\Xi,\boldsymbol{\Xi},\mu;\mathbb{R})\times L^2(\Xi,\boldsymbol{\Xi},\mu;\mathbb{R})\to\mathbb{R}\qquad:\Leftrightarrow\qquad\langle f,g\rangle=\int fg\,\mathrm{d}\mu.
\end{equation*}
This \emph{random function space} (aka RFS in this paper) is the space of all measurable functions $f:(\Xi,\boldsymbol{\Xi})\to(\mathbb{R},\mathcal{B}_\mathbb{R})$ that are square-integrable with respect to $\mu$.
(By $f$ we actually mean an equivalence class of square-integrable functions that are equal $\mu$-almost everywhere; usually denoted by $[f]$ in the literature.) This inner product space is known to form a Hilbert space because it is complete under the metric induced by the inner product.
Furthermore, we define $\{\Psi_j:(\Xi,\boldsymbol{\Xi})\to(\mathbb{R},\mathcal{B}_\mathbb{R})\}_{j=0}^\infty$ to be a complete orthogonal basis in $\mathscr{Z}$, such that $\Psi_0(\xi)=1$ for all $\xi\in\Xi$.
It is worth noting that such a basis does not necessarily need to consist of $d$-variate polynomials as in Refs.~\cite{xiu2002wiener,xiu2010numerical}, but may also consist of more general functions (including non-elementary functions such as wavelets).
\end{definition}

Therefore, any function $f\in\mathscr{Z}$ can be represented in a Fourier series of the form:
\begin{equation*}
f=\sum_{j=0}^\infty f^j\Psi_j,
\end{equation*}
where $f^j$ denotes the $j$-th coefficient of the series with the superscript not symbolizing an exponentiation.

Moreover, the dual space of $\mathscr{Z}$, which we denote by $\mathscr{Z}'$, is simply the space spanned by the continuous linear functionals $\{\Psi^i:\mathscr{Z}\to\mathbb{R}\}_{i=0}^\infty$ defined by:
\begin{equation*}
\Psi^i[f]:=[\Psi^i,f]=\frac{\langle\Psi_i,f\rangle}{\langle\Psi_i,\Psi_i\rangle}\equiv f^i,
\end{equation*}
where $[\,\cdot\,,\cdot\,]:\mathscr{Z}'\times\mathscr{Z}\to\mathbb{R}$ is the natural pairing map between $\mathscr{Z}$ and $\mathscr{Z}'$.
This continuous dual space is also known to form a Hilbert space, thanks to the Riesz representation theorem \cite{rudin1987real}.

We recall that the orthogonality property of the basis $\{\Psi_j\in\mathscr{Z}\}_{j=0}^\infty$ implies that:
\begin{equation*}
\langle\Psi_i,\Psi_j\rangle:=\int\Psi_i\Psi_j\,\mathrm{d}\mu=\langle\Psi_i,\Psi_i\rangle\,\delta_{ij},
\end{equation*}
where $\delta_{ij}$ is the Kronecker delta.

\begin{definition}[Solution space and root space]
Let $\mathscr{U}=\mathscr{T}(2)\otimes\mathscr{Z}$ and $\mathscr{V}=\mathscr{T}(0)\otimes\mathscr{Z}$ be, respectively, the \emph{solution space} and the \emph{root space} of the system.
\end{definition}

\begin{figure}
\centering
\includegraphics[]{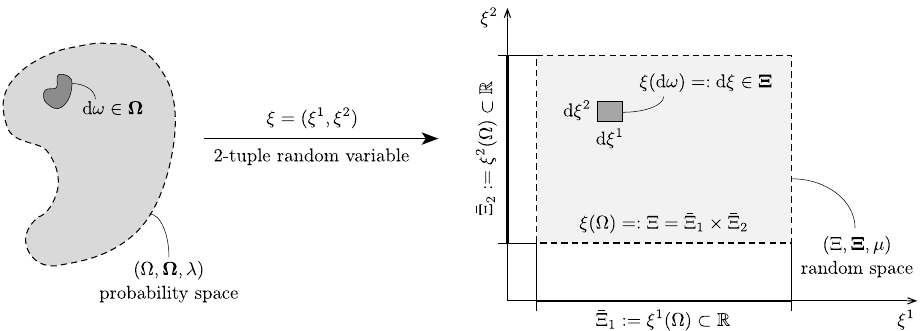}
\caption{Relationship between probability space and random space for $d=2$}
\label{fig1ProbabilityRandomSpace}
\end{figure}

In what follows, we assume that the components of the $d$-tuple random variable $\xi=(\xi^1,\ldots,\xi^d)$ are mutually independent, and as sketched in Fig.~\ref{fig1ProbabilityRandomSpace} that the random domain $\Xi$ is a hypercube of $d$ dimensions obtained by performing a $d$-fold Cartesian product of intervals $\bar{\Xi}_i:=\xi^i(\Omega)$.
It is for this reason that we define the probability measure in $\mathscr{Z}$ hereafter as
\begin{equation*}
\mu=\bigotimes_{i=1}^d\mu^i,\quad\text{or equivalently,}\quad\mathrm{d}\mu\equiv\mu(\mathrm{d}\xi)=\prod_{i=1}^d\mu^i(\mathrm{d}\xi^i)\equiv\mathrm{d}\mu^1\cdots\mathrm{d}\mu^d,
\end{equation*}
where $\mu^i(\mathrm{d}\xi^i)=:\mathrm{d}\mu^i$ represents the probability measure of $\mathrm{d}\xi^i$ in the vicinity of $\xi^i\in\bar{\Xi}_i$.

\section{Problem statement}\label{sec1ProSta}

In this work, we are interested in solving the non-autonomous, second-order ODE described below.

Find the real-valued stochastic process $u:\mathfrak{T}\times\Xi\to\mathbb{R}$ in $\mathscr{U}$, such that ($\mu$-a.e.):
\begin{subequations}\label{eq1ProSta1000}
\begin{align}
m\ddot{u}+\mathcal{F}[u,\dot{u}]=p&\qquad\text{on $\mathfrak{T}\times\Xi$}\label{eq1ProSta1000a}\\
\big\{u(0,\cdot\,)=\mathscr{u},\,\dot{u}(0,\cdot\,)=\mathscr{v}\big\}&\qquad\text{on $\{0\}\times\Xi$},\label{eq1ProSta1000b}
\end{align}
\end{subequations}
where $m,\mathcal{F}[u,\dot{u}],p:\mathfrak{T}\times\Xi\to\mathbb{R}$ are elements of $\mathscr{V}$ such that $m(t,\xi)\neq0$ for all $(t,\xi)\in\mathfrak{T}\times\Xi$, and $\mathscr{u},\mathscr{v}:\Xi\to\mathbb{R}$ are elements of $\mathscr{Z}$.
Note that in \eqref{eq1ProSta1000}, $\dot{u}:=\partial_t u$ and $\ddot{u}:=\partial_t^2 u$ denote, respectively, the first and second partial derivatives of $u$ with respect to time.
Thus, $\dot{u}:\mathfrak{T}\times\Xi\to\mathbb{R}$ is an element of $\mathscr{T}(1)\otimes\mathscr{Z}$, and $\ddot{u}:\mathfrak{T}\times\Xi\to\mathbb{R}$ is an element of $\mathscr{T}(0)\otimes\mathscr{Z}\equiv\mathscr{V}$ (the root space).

When \eqref{eq1ProSta1000} is specialized to study the (nonlinear) behavior of a single-degree-of-freedom system, $m\ddot{u}$ represents the inertial force of the system with $m:\mathfrak{T}\times\Xi\to\mathbb{R}^+$ symbolizing the mass of the system, $\mathcal{F}[u,\dot{u}]$ is the damping and resisting force of the system, and $p$ is the external force acting on the system.
Furthermore, in this case $u$, $\dot{u}$ and $\ddot{u}$ denote the displacement, the velocity and the acceleration response of the system, respectively.

When \eqref{eq1ProSta1000} is written in modeling notation, it becomes
\begin{equation}\label{eq1ProSta1000star}
y=\boldsymbol{\mathcal{M}}[u][x]\quad\text{subject to initial condition}\quad\boldsymbol{\mathcal{I}}[u],\tag{\ref{eq1ProSta1000}*}
\end{equation}
where $\boldsymbol{\mathcal{M}}[u]:\mathscr{V}^3\to\mathscr{V}^s$ represents the mathematical model of the system defined by \eqref{eq1ProSta1000a}, $x=(x_1,x_2,x_3):\mathfrak{T}\times\Xi\to\mathbb{R}^3$ is the $3$-tuple input of $\boldsymbol{\mathcal{M}}[u]$, and $y=(y_1,\ldots,y_s):\mathfrak{T}\times\Xi\to\mathbb{R}^s$ is the $s$-tuple output of $\boldsymbol{\mathcal{M}}[u]$ (aka the $s$-tuple observable in physics or the $s$-tuple response in engineering).
In addition, $\boldsymbol{\mathcal{I}}[u]$ represents the initial condition for $\boldsymbol{\mathcal{M}}[u]$ which is given by \eqref{eq1ProSta1000b}.
Thus, by comparing \eqref{eq1ProSta1000} to \eqref{eq1ProSta1000star}, the components of $x$ can be identified in the following way: $x_1=m$, $x_2=\mathcal{F}[u,\dot{u}]$ and $x_3=p$.
The objective of this mathematical model is to propagate and quantify the effects of input uncertainty $x$ on system's output $y$.
Therefore, besides seeking $u$ in $\mathscr{U}$ as mentioned earlier, it is also important to compute the probability moments of $y$ as time progresses.

\section{Solution based on the spectral approach}\label{sec1SolSpeApp}

\subsection{Discretization of random function space}

Since by hypothesis $u$ is an element of $\mathscr{U}$, then it can be represented by the Fourier series
\begin{equation}\label{eq1SolSpeApp1000}
u(t,\xi)=\sum_{j=0}^\infty u^j(t)\,\Psi_j(\xi),
\end{equation}
where $u^j$ is a temporal function in $\mathscr{T}(2)$ denoting the $j$-th random mode of $u$.
This series, usually referred to as \emph{stochastic spectral expansion} in the literature \cite{le2010spectral,sullivan2015introduction}, has the remarkable property that when $u$ is sufficiently smooth in the solution space (and, of course, provided that the basis is orthogonal with respect to the probability measure defined in $\mathscr{Z}$), it leads to exponential convergence to the solution \cite{xiu2002wiener,ernst2012convergence}.

For the purpose of this manuscript, let us simply consider a $p$-discretization of the random function space as follows.
Let $\mathscr{Z}^{[P]}=\mathrm{span}\{\Psi_j\}_{j=0}^P$ be a finite subspace of $\mathscr{Z}$ with $P+1\in\mathbb{N}_1$ denoting the dimensionality of the subspace.
If we let $u^{[P]}(t,\cdot\,)$ be an element of $\mathscr{Z}^{[P]}$, then it is evident from \eqref{eq1SolSpeApp1000} that\footnote{As long as we assume that $\{\Psi_j\}_{j=0}^\infty$ is well-graded to carry out the approximation of $u$ this way.}:
\begin{equation}\label{eq1SolSpeApp1030}
u(t,\xi)\approx u^{[P]}(t,\xi)=\sum_{j=0}^P u^j(t)\,\Psi_j(\xi)\equiv u^j(t)\,\Psi_j(\xi),
\end{equation}
where for notational convenience we have omitted the summation sign in the last equality (aka Einstein summation convention).
Therefore, unless otherwise noted hereinafter, a summation sign will always be implied over the repeated index $j\in\{0,1,\ldots,P\}$.

Substituting \eqref{eq1SolSpeApp1030} into \eqref{eq1ProSta1000} gives
\begin{subequations}\label{eq1SolSpeApp1040}
\begin{align}
m\ddot{u}^j\Psi_j+\mathcal{F}[u^j\Psi_j,\dot{u}^j\Psi_j]=p&\qquad\text{on $\mathfrak{T}\times\Xi$}\label{eq1SolSpeApp1040a}\\
\big\{u^j(0)\,\Psi_j=\mathscr{u},\,\,\dot{u}^j(0)\,\Psi_j=\mathscr{v}\big\}&\qquad\text{on $\{0\}\times\Xi$}.\label{eq1SolSpeApp1040b}
\end{align}
\end{subequations}

Projecting \eqref{eq1SolSpeApp1040} onto $\mathscr{Z}^{[P]}$ yields a system of $P+1$ ordinary differential equations of second order in the variable $t$, where the unknowns are the random modes $u^j=u^j(t)$ and $\dot{u}^j=\dot{u}^j(t)$:
\begin{subequations}\label{eq1SolSpeApp1050}
\begin{align}
\langle\Psi_i,m\Psi_j\rangle\,\ddot{u}^j+\langle\Psi_i,\mathcal{F}[u^j\Psi_j,\dot{u}^j\Psi_j]\rangle=\langle\Psi_i,p\rangle &\qquad\text{on $\mathfrak{T}$}\label{eq1SolSpeApp1050a}\\
\big\{u^i(0)=\langle\Psi_i,\mathscr{u}\rangle/\langle\Psi_i,\Psi_i\rangle,\,\,\dot{u}^i(0)=\langle\Psi_i,\mathscr{v}\rangle/\langle\Psi_i,\Psi_i\rangle\big\}&\qquad\text{on $\{0\}$}\label{eq1SolSpeApp1050b}
\end{align}
\end{subequations}
with $i,j\in\{0,1,\ldots,P\}$.
Note that in order to get \eqref{eq1SolSpeApp1050}, we simply applied on both sides of each equation the linear functionals $\{\Psi^i\in\mathscr{Z}'\}_{i=0}^P$ one by one, and then we simplified the resulting expressions.
It is also worth noting that because the randomness of the stochastic system has effectively been absorbed by the application of the aforementioned functionals, the system of equations that we are dealing with at this point is no longer `stochastic' but `deterministic'.
In other words, the system now depends merely on the time variable $t$ rather than on the tuple $(t,\xi)$.

System \eqref{eq1SolSpeApp1050} can also be restated using multilinear and tensor algebra notation as follows:
\begin{subequations}\label{eq1SolSpeApp1060}
\begin{align}
m\indices{^i_j}\ddot{u}^j+\mathcal{F}^i[u^j,\dot{u}^j]=p^i &\qquad\text{on $\mathfrak{T}$}\label{eq1SolSpeApp1060a}\\
\big\{u^i(0)=\mathscr{u}^i,\,\dot{u}^i(0)=\mathscr{v}^i\big\}&\qquad\text{on $\{0\}$},\label{eq1SolSpeApp1060b}
\end{align}
\end{subequations}
where $i,j\in\{0,1,\ldots,P\}$ (summation sign implied over repeated index $j$), and:
\begin{gather*}
m\indices{^i_j}(t)=\langle\Psi_i,m(t,\cdot\,)\,\Psi_j\rangle/\langle\Psi_i,\Psi_i\rangle,\quad
\mathcal{F}^i[u^j,\dot{u}^j](t)=\langle\Psi_i,\mathcal{F}[u^j\Psi_j,\dot{u}^j\Psi_j](t,\cdot\,)\rangle/\langle\Psi_i,\Psi_i\rangle,\\
p^i(t)=\langle\Psi_i,p(t,\cdot\,)\rangle/\langle\Psi_i,\Psi_i\rangle,\quad
\mathscr{u}^i=\langle\Psi_i,\mathscr{u}\rangle/\langle\Psi_i,\Psi_i\rangle\quad\text{and}\quad
\mathscr{v}^i=\langle\Psi_i,\mathscr{v}\rangle/\langle\Psi_i,\Psi_i\rangle,
\end{gather*}
whence $m\indices{^i_j},\mathcal{F}^i[u^j,\dot{u}^j],p^i\in\mathscr{T}(0)$ and $\mathscr{u}^i,\mathscr{v}^i\in\mathbb{R}$.
To simplify notation, we have taken $\mathcal{F}^i[u^j,\dot{u}^j]$ as the short notation for $\mathcal{F}^i[u^0,\ldots,u^j,\ldots,u^P,\dot{u}^0,\ldots,\dot{u}^j,\ldots,\dot{u}^P]$.

To evaluate the inner products approximately, any integration technique of one's choice can be used, including those addressed in \cite{davis2007methods,novak1996high,bungartz2004sparse,leobacher2014introduction,guo2020constructing}.
If, for instance, a Gaussian quadrature rule is adopted, the inner products are computed with:
\begin{equation*}
\langle f,g\rangle:=\int fg\,\mathrm{d}\mu\approx\mathcal{Q}^{[Q]}[fg]:=\sum_{i=1}^Q f(\xi_i)\,g(\xi_i)\, w_i,
\end{equation*}
where $w_i\in\mathbb{R}^+$ denotes the quadrature weight associated with the Gaussian quadrature point $\xi_i\in\Xi$ (w.r.t.~$\mu$), and $Q\in\mathbb{N}_1$ represents the number of quadrature points involved in approximating the evaluation of the inner product.

\subsection{Discretization of temporal function space}\label{sec1DisTemFunSpa}

Once $\mathscr{Z}$ has been discretized, the temporal function space can be discretized using an $(h,p)$-discretization for $\mathscr{T}(2)$.
In the literature, there exists an extensive list of time integration techniques that one can employ in order to solve the ODE system given by \eqref{eq1SolSpeApp1060} numerically. For example, the Runge-Kutta method \cite{butcher1987numerical} of fourth-order (aka RK4 method) or the Newmark-$\beta$ method \cite{newmark1959method,chopra2012dynamics}.

\subsection{Computation of probability moments}

In this manuscript, the probability moments of interest are the mean and the variance of the system's response.
For this reason, we define these objects below.

Suppose that $z:=y_k$ is the $k$-th component of output $y=\boldsymbol{\mathcal{M}}[u][x]$.
If $z\in\mathscr{V}$, then it can approximately be expanded with a truncated Fourier series similar to the one set forth in \eqref{eq1SolSpeApp1030} to obtain:
\begin{equation}\label{eq1SolSpeApp1270}
z(t,\xi)\approx z^{[P]}(t,\xi)=\sum_{j=0}^P z^j(t)\,\Psi_j(\xi)\equiv z^j(t)\,\Psi_j(\xi),
\end{equation}
where the $j$-th random mode of $z$ is given by
\begin{equation*}
z^j(t)=\frac{\langle\Psi_j,z(t,\cdot\,)\rangle}{\langle\Psi_j,\Psi_j\rangle}.
\end{equation*}
(Note that $P$ in expression \eqref{eq1SolSpeApp1270} does not necessarily need to be the same as in \eqref{eq1SolSpeApp1030}.)

The expectation of $z$, $\mathbf{E}[z]:\mathfrak{T}\to\mathbb{R}$, is simply given by the first random mode of $z$:
\begin{equation}\label{eq1SolSpeApp1280}
\mathbf{E}[z](t):=\int z(t,\cdot\,)\,\mathrm{d}\mu=z^0(t),
\end{equation}
whereas the variance of $z$, $\mathrm{Var}[z]:\mathfrak{T}\to\mathbb{R}^+_0$, is defined by the partial sum:
\begin{equation}\label{eq1SolSpeApp1290}
\mathrm{Var}[z](t):=\int (z(t,\cdot\,)-\mathbf{E}[z](t))^2\,\mathrm{d}\mu=\sum_{j=1}^P\langle\Psi_j,\Psi_j\rangle\,z^j(t)\,z^j(t).
\end{equation}

\section{Flow-driven spectral chaos (FSC) method}\label{sec1FSCMet}

\subsection{Stochastic flow map}\label{subsec1StoFloMap}

Observe that the stochastic system given by \eqref{eq1ProSta1000} can also be expressed as:
\begin{subequations}\label{eq1FSCMet1000}
\begin{align}
\partial_t^2 u(t,\xi):=\mathscr{f}(t,\xi,s(t,\xi))=\big(p(t,\xi)-\mathcal{F}[u,\dot{u}](t,\xi)\big)/m(t,\xi) &\qquad\text{on $\mathfrak{T}\times\Xi$}\label{eq1FSCMet1000a}\\
\big\{u(0,\xi)=\mathscr{u}(\xi),\,\dot{u}(0,\xi)=\mathscr{v}(\xi)\big\} &\qquad\text{on $\{0\}\times\Xi$},\label{eq1FSCMet1000b}
\end{align}
\end{subequations}
where $s=(u,\dot{u})\in\prod_{j=1}^2\mathscr{T}(3-j)\otimes\mathscr{Z}$ is the configuration state of the system over $\mathfrak{T}\times\Xi$, and $\mathscr{f}:\mathfrak{T}\times\Xi\times\mathbb{R}^2\to\mathbb{R}$ is a noisy, non-autonomous function defining the response $\ddot{u}=\partial_t^2 u$.

Therefore, if the solution is analytic on $\mathfrak{T}$ for all $\xi\in\Xi$, then it can be represented by the Taylor series:
\begin{equation}\label{eq1FSCMet1010}
u(t_i+h,\xi)=\sum_{j=0}^\infty \frac{h^j}{j!}\partial_t^j u(t_i,\xi)=\sum_{j=0}^M\frac{h^j}{j!}\partial_t^j u(t_i,\xi)+O(h^{M+1})(\xi),
\end{equation}
where $h:=t-t_i$ is the time-step size used for the simulation around $t_i$ (once $t$ is fixed), $t_i\in\mathfrak{T}$ is the time instant of the simulation, and $M\in\mathbb{N}_1$ is the order of the flow map we are interested to implement.

For this system, the \emph{stochastic flow map} of order $M$, $\varphi(M):\mathbb{R}\times\mathscr{Z}^2\to\mathscr{Z}^2$, can be defined as a random map given by:
\begin{multline}\label{eq1FSCMet1040}
\varphi(M)(h,s(t_i,\cdot\,))=:s(t_i+h,\cdot\,)=\big(u(t_i+h,\cdot\,),\dot{u}(t_i+h,\cdot\,)\big)-O(h^{M+1})\\
=\bigg(\sum_{j=0}^M\frac{h^j}{j!}\partial_t^j u(t_i,\cdot\,),\sum_{j=0}^M\frac{h^j}{j!}\partial_t^{j+1} u(t_i,\cdot\,)\bigg),
\end{multline}
where this $s(t_i+h,\cdot\,)$ is the same as in \eqref{eq1FSCMet1000} if $M\to\infty$. However, to avoid unnecessary complexity in notation, no distinction between these definitions will be made in this work. That is, from now on we will assume that there is an equivalence relation $\sim$ between $s(t_i+h,\cdot\,)-O(h^{M+1})$ and the $s(t_i+h,\cdot\,)$ defined in \eqref{eq1FSCMet1040}. Notice that \eqref{eq1FSCMet1040} now requires that $s=(u,\dot{u})\in\prod_{j=1}^2\mathscr{T}(M-j+2)\otimes\mathscr{Z}$.

For the sake of illustration, suppose $M=4$.
Differentiating \eqref{eq1FSCMet1000a} with respect to time three times (i.e.~$M-1$ times) gives
\begin{subequations}\label{eq1FSCMet1045}
\begin{gather}
\partial_t^3 u:=\mathrm{D}_t\mathscr{f}=\partial_t\mathscr{f}+\partial_u\mathscr{f}\,\partial_t u+\partial_{\dot{u}}\mathscr{f}\,\partial_t^2 u\\
\partial_t^4 u:=\mathrm{D}_t^2\mathscr{f}=\partial_t^2\mathscr{f}+2\,\partial_{tu}^2\mathscr{f}\,\partial_t u+(2\,\partial_{t\dot{u}}^2\mathscr{f}+\partial_u\mathscr{f})\,\partial_t^2 u+\partial_{\dot{u}}\mathscr{f}\,\partial_t^3 u\\
\partial_t^5 u:=\mathrm{D}_t^3\mathscr{f}=\partial_t^3\mathscr{f}+3\,\partial_{ttu}^3\mathscr{f}\,\partial_t u+3\,(\partial_{tt\dot{u}}^3\mathscr{f}+\partial_{tu}^2\mathscr{f})\,\partial_t^2 u+(3\,\partial_{t\dot{u}}^2\mathscr{f}+\partial_u\mathscr{f})\,\partial_t^3 u+\partial_{\dot{u}}\mathscr{f}\,\partial_t^4 u.
\end{gather}
\end{subequations}
Hence, when $M=4$, the stochastic flow map of the system, $\varphi(4)=(\varphi^1(4),\varphi^2(4))$, is given by:
\begin{multline}\label{eq1FSCMet1050}
\varphi(4)(h,s(t_i,\cdot\,))=:s(t_i+h,\cdot\,)=\big(u(t_i+h,\cdot\,),\dot{u}(t_i+h,\cdot\,)\big)-O(h^5)\\
=\bigg(\sum_{j=0}^4\frac{h^j}{j!}\partial_t^j u(t_i,\cdot\,),\sum_{j=0}^4\frac{h^j}{j!}\partial_t^{j+1} u(t_i,\cdot\,)\bigg),
\end{multline}
where the second and higher time derivatives of $u$ at $t=t_i$ are computed with \eqref{eq1FSCMet1000a} and \eqref{eq1FSCMet1045}, respectively.

Note that if \eqref{eq1FSCMet1000} is an autonomous ODE, the expressions prescribed by \eqref{eq1FSCMet1045} reduce to:
\begin{equation}
\partial_t^3 u=\partial_u\mathscr{f}\,\partial_t u+\partial_{\dot{u}}\mathscr{f}\,\partial_t^2 u,\quad\partial_t^4 u=\partial_u\mathscr{f}\,\partial_t^2u+\partial_{\dot{u}}\mathscr{f}\,\partial_t^3u\quad\text{and}\quad\partial_t^5u=\partial_u\mathscr{f}\,\partial_t^3u+\partial_{\dot{u}}\mathscr{f}\,\partial_t^4u.\tag{\ref{eq1FSCMet1045}*}	
\end{equation}

\subsection{Enriched stochastic flow map}

In this work, we define the \emph{enriched stochastic flow map} of order $M$, $\hat{\varphi}(M):\mathbb{R}\times\mathscr{Z}^{M+2}\to\mathscr{Z}^{M+2}$, such that its $k$-th component is given by:
\begin{equation}
\hat{\varphi}^k(M)(h,\hat{s}(t_i,\cdot\,))=:\hat{s}^k(t_i+h,\cdot\,)=%
\begin{cases}
\varphi^k(M)(h,s(t_i,\cdot\,))& \text{for $k\in\{1,2\}$}\\
\mathrm{D}_t^{k-3}\mathscr{f}(t_i+h,\cdot\,,s(t_i+h,\cdot\,)) & \text{otherwise,}
\end{cases}
\end{equation}
where $\hat{s}=(u,\dot{u},\ldots,\partial_t^{M+1}u)\in\prod_{j=1}^{M+2}\mathscr{T}(M-j+2)\otimes\mathscr{Z}$ is called the \emph{enriched configuration state} of the system over $\mathfrak{T}\times\Xi$, and $k\in\{1,2,\ldots,M+2\}$.

For instance, when $M=4$, the components of the enriched stochastic flow map, $\hat{\varphi}(4)$, are
\begin{subequations}
\begin{gather}
\hat{\varphi}^1(4)(h,s(t_i,\cdot\,)):=\varphi^1(4)(h,s(t_i,\cdot\,))=s^1(t_i+h,\cdot\,)=u(t_i+h,\cdot\,)-O(h^5),\\
\hat{\varphi}^2(4)(h,s(t_i,\cdot\,)):=\varphi^2(4)(h,s(t_i,\cdot\,))=s^2(t_i+h,\cdot\,)=\dot{u}(t_i+h,\cdot\,)-O(h^5),\\
\hat{\varphi}^3(4):=\mathscr{f}=\partial_t^2 u,\quad \hat{\varphi}^4(4):=\mathrm{D}_t\mathscr{f}=\partial_t^3 u,\quad \hat{\varphi}^5(4):=\mathrm{D}_t^2\mathscr{f}=\partial_t^4 u\quad\text{and}\quad\hat{\varphi}^6(4):=\mathrm{D}_t^3\mathscr{f}=\partial_t^5 u,
\end{gather}
\end{subequations}
where $\hat{\varphi}^1(4)$ and $\hat{\varphi}^2(4)$ are computed with \eqref{eq1FSCMet1050}, $\hat{\varphi}^3(4)$ with \eqref{eq1FSCMet1000a}, and $\{\hat{\varphi}^k(4)\}_{k=4}^6$ with \eqref{eq1FSCMet1045}.

\subsection{Derivation of the FSC method}

According to Section \ref{subsec1StoFloMap}, the state of a system driven by a stochastic flow map of order $M$ is:
\begin{equation}\label{eq1FSCMet2010}
u(t,\xi)=\sum_{j=0}^M\frac{(t-t_i)^j}{j!}\partial_t^j u(t_i,\xi)\quad\text{and}\quad \dot{u}(t,\xi)=\sum_{j=0}^M\frac{(t-t_i)^j}{j!}\partial_t^{j+1} u(t_i,\xi),
\end{equation}
with the provision that the stochastic process $u$ is analytic on the temporal domain.

From these two expressions, it can be seen that the state of the system has been decomposed effectively into deterministic and non-deterministic parts.
That is, the deterministic part $(t-t_i)^j/j!$ consisting of a temporal function in $\mathscr{T}$, and the non-deterministic part $\partial_t^j u(t_i,\xi)$ consisting of a random function in $\mathscr{Z}$.
If the set of functions associated with the non-deterministic part, i.e.~$\{\partial_t^j u(t_i,\cdot\,)\}_{j=0}^{M+1}$, is orthogonalized with respect to the measure in $\mathscr{Z}$, then \eqref{eq1FSCMet2010} can also be written in the following way:
\begin{equation}\label{eq1FSCMet2020}
u(t,\xi)=\sum_{j=1}^{M+2} u^j(t)\,\Psi_j(\xi)\quad\text{and}\quad\dot{u}(t,\xi)=\sum_{j=1}^{M+2} \dot{u}^j(t)\,\Psi_j(\xi).
\end{equation}
Hence, if the space associated with the stochastic part of the solution space were to be spanned with $\{\Psi_j\}_{j=1}^{M+2}$, then $u(t,\cdot\,)$ and $\dot{u}(t,\cdot\,)$ would be elements of that space around $t=t_i$.
However, since one cannot always guarantee that the constant functions are in such a space construction, we prefer to write \eqref{eq1FSCMet2020} in the following final form instead:
\begin{equation}\label{eq1FSCMet2030}
u(t,\xi)=\sum_{j=0}^{M+2} u^j(t)\,\Psi_j(\xi)\quad\text{and}\quad\dot{u}(t,\xi)=\sum_{j=0}^{M+2} \dot{u}^j(t)\,\Psi_j(\xi),
\end{equation}
where $\Psi_0\equiv1$ is the identically-equal-to-one function as per Definition \ref{sec1defRanFunSpa}.

Therefore, for a system driven by a stochastic flow map of order $M$, the maximum number of basis vectors to use in a simulation with FSC is bounded from above by $M+3$.
Hence, regardless of the dimensionality of the random space, the probability information of the system's state can be completely captured in $\mathscr{Z}^{[P]}$ if $P=M+2$.
It is for this reason that our FSC scheme is superior in terms of efficiency in comparison to mTD-gPC which uses a combination of full and total-order tensor products to construct a suitable basis for $\mathscr{Z}^{[P]}$ around $t=t_i$.
We emphasize, however, that the FSC scheme does not address by itself the curse of dimensionality at the random-space level, since we still have the issue that the bigger the random space is ($d\gg1$), the more difficult it is to compute the inner products accurately in \eqref{eq1SolSpeApp1060} and \eqref{eq1SolSpeApp1290}.
This is still an open area of research and there are several approaches available in the literature for dealing with this issue \cite{davis2007methods,novak1996high,bungartz2004sparse,leobacher2014introduction,guo2020constructing}.

Moreover, to reduce the computational cost associated with orthogonalizing $M+2$ basis vectors, it is sometimes convenient to start the FSC analysis with the smallest value for $M$ (i.e.~$M=1$), and then to progressively increment it if more accurate results are desired for the simulation.
Therefore, the minimum number of basis vectors to use in a simulation with FSC is bounded from below by 4.

\subsection{FSC scheme}\label{sec1FSCSch}

\begin{figure}
\centering
\includegraphics[]{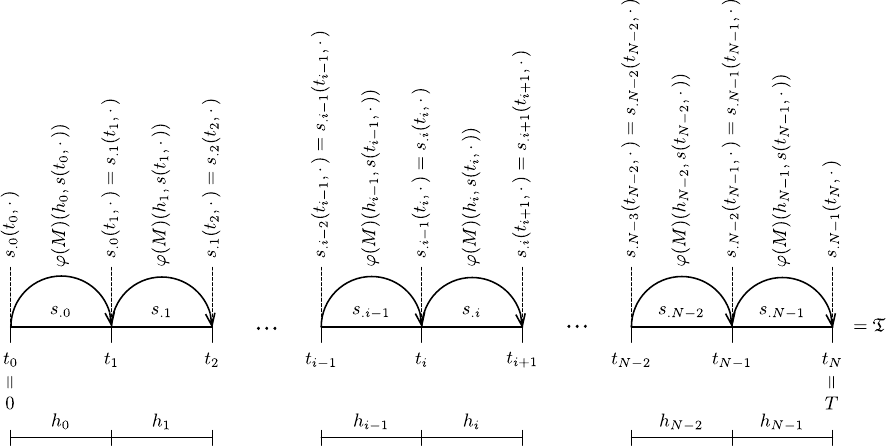}
\caption{Evolution of a dynamical system via a stochastic flow map of order $M$, provided that $h_i$ is taken sufficiently small and $M$ sufficiently large so as to have $s(t_N,\cdot\,)\approx s_{.N-1}(t_N,\cdot\,)$ at the end of the simulation. Then, as a means to avoid complexity in notation, we take $s(t_i+h,\cdot\,)\sim s_{.i}(t_i+h,\cdot\,)$ for $h\geq0$, and $s(t_i+h,\cdot\,)\sim s_{.i-1}(t_i+h,\cdot\,)$ for $h\leq0$}
\label{fig1StochasticFlowMap}
\end{figure}

Suppose that a stochastic system such as \eqref{eq1ProSta1000} has been given.
Let $\{\mathfrak{T}_i\}_{i=0}^{N-1}$ be a partition of the temporal domain, where $\mathfrak{T}_i\neq\emptyset$ represents the $i$-th interval of the partition, and define
$s_{.i}=s|_{\mathrm{cl}(\mathfrak{T}_i)\times\Xi}$ to be the restriction of $s$ to $\mathfrak{R}_i:=\mathrm{cl}(\mathfrak{T}_i)\times\Xi$.
(Recall that $s=(u,\dot{u})$ represents the configuration state of the system over $\mathfrak{T}\times\Xi$.)
Then, if the system is driven by a stochastic flow map of order $M$ (Fig.~\ref{fig1StochasticFlowMap}), proceed as follows:
\begin{enumerate}
\item Loop across the temporal domain from $i=0$ to $i=N-1$.
\begin{enumerate}
\item Define a solution representation for the configuration state $s_{.i}$ in the following way.
\begin{itemize}
\item Take $\{\Phi_{j.i}:=\hat{\varphi}^j(M)(0,\hat{s}(t_i,\cdot\,))\}_{j=1}^P$ to be an ordered set of linearly independent functions in $\mathscr{Z}$ with $3\leq P\leq M+2$, and define $\Phi_{0.i}\equiv 1$.
Observe that $\hat{\varphi}(M)(0,\hat{s}(t_i,\cdot\,))\equiv\hat{s}_{.i}(t_i,\cdot\,)=\hat{s}_{.i-1}(t_i,\cdot\,)$ for $i\geq1$.
However, if $i=0$, then $\hat{\varphi}(M)(0,\hat{s}(t_0,\cdot\,))\equiv \hat{s}(0,\cdot\,)$.
(\emph{Note:} When the initial conditions are deterministic or linearly dependent, please see Remark \ref{rmk1paper100}.)
\item Orthogonalize the set $\{\Phi_{j.i}\}_{j=0}^P$ using the Gram-Schmidt process \cite{cheney2010linear}, so that the resulting set $\{\Psi_{j.i}\}_{j=0}^P$ is an orthogonal basis in $\mathscr{Z}$.
That is, for $j\in\{0,1,\ldots,P\}$:
\begin{equation}\label{eq1FSCMet3010}
\Psi_{j.i}:=\Phi_{j.i}-\sum_{k=0}^{j-1}\frac{\langle\Phi_{j.i},\Psi_{k.i}\rangle}{\langle\Psi_{k.i},\Psi_{k.i}\rangle}\Psi_{k.i}.
\end{equation}
\item Define $\mathscr{Z}_i^{[P]}=\mathrm{span}\{\Psi_{j.i}\}_{j=0}^P$ to be a $p$-discretization of $\mathscr{Z}$ over the region $\mathfrak{R}_i$.
Since $\mathscr{Z}_i^{[P]}$ is an evolving function space, expansion \eqref{eq1SolSpeApp1030} is now to be read as:
\begin{equation}\label{eq1FSCMet3020}
u_{.i}(t,\xi)\approx u^{[P]}_{.i}(t,\xi)=\sum_{j=0}^P u{^j}_{\!.i}(t)\,\Psi_{j.i}(\xi)\equiv u{^j}_{\!.i}(t)\,\Psi_{j.i}(\xi).\tag{\ref{eq1SolSpeApp1030}*}
\end{equation}
Hence, $\dot{u}_{.i}=\partial_t u_{.i}$.
\end{itemize}
\item Transfer the random modes of $s_{.i-1}=(u_{.i-1},\dot{u}_{.i-1})$ to $s_{.i}=(u_{.i},\dot{u}_{.i})$ at $t=t_i$, given that $i\geq1$.\par
One way to achieve this is to ensure that the probability information of the system's state is transferred in the mean-square sense.
Put differently, we wish to make sure that the equalities shown below hold in the mean-square sense (summation sign implied only over repeated index $k$):
\begin{subequations}\label{eq1FSCMet3030}
\begin{gather}
u_{.i}(t_i,\xi)=u_{.i-1}(t_i,\xi)\quad\iff\quad u{^k}_{\!.i}(t_i)\,\Psi_{k.i}(\xi)=u{^k}_{\!.i-1}(t_i)\,\Psi_{k.i-1}(\xi)\,\,\\
\dot{u}_{.i}(t_i,\xi)=\dot{u}_{.i-1}(t_i,\xi)\quad\iff\quad \dot{u}{^k}_{\!.i}(t_i)\,\Psi_{k.i}(\xi)=\dot{u}{^k}_{\!.i-1}(t_i)\,\Psi_{k.i-1}(\xi).
\end{gather}
\end{subequations}
Projecting \eqref{eq1FSCMet3030} onto $\mathscr{Z}_i^{[P]}$ gives:
\begin{subequations}
\begin{gather}
u{^k}_{\!.i}(t_i)\,\langle\Psi_{j.i},\Psi_{k.i}\rangle=u{^k}_{\!.i-1}(t_i)\,\langle\Psi_{j.i},\Psi_{k.i-1}\rangle\,\,\\
\dot{u}{^k}_{\!.i}(t_i)\,\langle\Psi_{j.i},\Psi_{k.i}\rangle=\dot{u}{^k}_{\!.i-1}(t_i)\,\langle\Psi_{j.i},\Psi_{k.i-1}\rangle.
\end{gather}
\end{subequations}
Thus, upon simplification yields the random modes of $s_{.i}=(u_{.i},\dot{u}_{.i})$ at $t=t_i$:
\begin{equation}\label{eq1FSCMet3040}
u{^j}_{\!.i}(t_i)=\frac{\langle\Psi_{j.i},\Psi_{k.i-1}\rangle}{\langle\Psi_{j.i},\Psi_{j.i}\rangle} u{^k}_{\!.i-1}(t_i)\quad\text{and}\quad 
\dot{u}{^j}_{\!.i}(t_i)=\frac{\langle\Psi_{j.i},\Psi_{k.i-1}\rangle}{\langle\Psi_{j.i},\Psi_{j.i}\rangle} \dot{u}{^k}_{\!.i-1}(t_i),\tag{\ref{eq1SolSpeApp1060b}*}
\end{equation}
where $j\in\{0,1,\ldots,P\}$ (summation sign implied over repeated index $k$).
These are to be interpreted as the initial conditions of the system over the region $\mathfrak{R}_i$.\par
If $i=0$, the initial conditions are computed with \eqref{eq1SolSpeApp1060b} directly.
\item Substitute \eqref{eq1FSCMet3020} into \eqref{eq1ProSta1000} to obtain \eqref{eq1SolSpeApp1040}.
\item Project \eqref{eq1SolSpeApp1040a} onto $\mathscr{Z}_i^{[P]}$ to obtain \eqref{eq1SolSpeApp1060a} subject to \eqref{eq1FSCMet3040}.
Note that if $i=0$, \eqref{eq1SolSpeApp1060a} is subject to \eqref{eq1SolSpeApp1060b}.
\item Integrate \eqref{eq1SolSpeApp1060} over time, as long as a suitable time integration method has been selected for solving the resulting system of equations. This step requires to find the random modes $\{u{^j}_{\!.i}(t_{i+1})\}_{j=0}^P$ and $\{\dot{u}{^j}_{\!.i}(t_{i+1})\}_{j=0}^P$.
\item Compute both the mean and the variance of each of the components of output $y=\boldsymbol{\mathcal{M}}[u][x]$ over $\mathfrak{R}_i$, by recurring to the formulas prescribed by \eqref{eq1SolSpeApp1280} and \eqref{eq1SolSpeApp1290}.
\end{enumerate}
\item Post-process results.
\end{enumerate}

\begin{remark}\label{rmk1paper100}
Any of the following two approaches can be carried out at the start of the simulation ($i=0$) to address the case when the initial conditions are deterministic or, more generally, linearly dependent:
\begin{itemize}
\item When the initial conditions are deterministic, the first two vectors are required to be removed from the set $\{\Phi_{j.i}:=\hat{\varphi}^j(M)(0,\hat{s}(t_i,\cdot\,))\}_{j=1}^P$ for they are constant, and when the initial conditions are stochastic but linearly dependent, only one of them needs to be removed from the set.
\item When the initial conditions are deterministic or linearly dependent, the gPC method \cite{xiu2002wiener,xiu2010numerical} can be employed instead to advance the state of the system one-time step forward; that is, from $s(t_0=0,\cdot\,)$ to $s(t_1,\cdot\,)$.
After this, the gPC method can be switched over FSC to continue pushing the system's state forward in time.
\end{itemize}
\end{remark}

\begin{remark}
Compared to the standard TD-gPC by Gerritsma et al.~\cite{gerritsma2010time}, in our FSC scheme we update the stochastic part of the solution space at every time step to minimize the error propagation over time.
We do this without loss of generality since the scheme can conveniently be modified to incorporate a stopping criterion of one's choice and update the random basis only when the criterion is met.
\end{remark}

\section{Numerical results}\label{sec1NumRes}

We demonstrate and compare the performance of the FSC scheme to the mTD-gPC scheme using two numerical examples for the dynamical system described in Section \ref{sec1ProSta}.
We also define the local error, $\epsilon:\mathscr{T}\to\mathscr{T}$, and the global error, $\epsilon_G:\mathscr{T}\to\mathbb{R}$, with the following expressions:
\begin{subequations}
\begin{gather}
\epsilon[f](t)=|f(t)-f_\mathrm{exact}(t)|\\
\epsilon_G[f]=\frac{1}{T}\int_\mathfrak{T} |f(t)-f_\mathrm{exact}(t)|\,\mathrm{d}t\approx\frac{\Delta t}{T}\sum_{i=0}^N|f(t_i)-f_\mathrm{exact}(t_i)|,
\end{gather}
\end{subequations}
where $\Delta t$ is the time-step size used for the simulation, $t_i\in\mathfrak{T}$ is the time instant of the simulation, and $N$ denotes the number of time steps employed in the simulation (with $t_0=0$ and $t_N=N\,\Delta t=T$).

\subsection{Single-degree-of-freedom system under free vibration}\label{sec1NumRes10}

We consider an undamped single-degree-of-freedom system governed by $m\ddot{u}+ku=0$ with mass $m=100$ kg and stochastic stiffness $k(\xi)=\xi$ subjected to free vibration.
Three different cases of stochasticity are considered as listed in Table \ref{tab1CaseStudies10}.
The system has an initial displacement of $u(0,\cdot\,)\equiv0.05$ m, and an initial velocity of $\dot{u}(0,\cdot\,)\equiv0.20$ m/s.
The duration of the simulation is set to $T=150$ s. To minimize the errors coming from the discretization of $\mathscr{T}(2)$, the time-step size is taken as $\Delta t=0.005$ s, meaning that a total of $N=30\,000$ time steps are employed in the simulation.
To integrate \eqref{eq1SolSpeApp1060} over time\footnote{In this problem we have taken $\mathcal{F}[u,\dot{u}]=ku$, and thus $\mathcal{F}^i[u^j,\dot{u}^j](t)=k\indices{^i_j}\,u^j(t)$ with $k\indices{^i_j}=\langle\Psi_i,k\Psi_j\rangle/\langle\Psi_i,\Psi_i\rangle$.}, we use the RK4 method, and as described in Section \ref{sec1FSCSch}, the stochastic part of the solution space is updated at every time step in order to obtain accurate results.
Lastly, because the initial conditions of the system are deterministic, we opt to employ the gPC method (with $P=6$) \cite{xiu2002wiener} for the first 5 seconds of the simulation, in an effort to ensure that the stochasticity of the system's state is well developed for the analysis with FSC or mTD-gPC.

\begin{table}
\centering\small
\caption{Case studies considered for a single-degree-of-freedom system under free vibration and stochastic stiffness}
\begin{tabular}{@{}p{0.125\textwidth}@{}p{0.35\textwidth}@{}p{0.25\textwidth}@{}p{0.275\textwidth}@{}}
\toprule
Case & Probability distribution$^*$ & Probability moments & Quadrature rule\\
\midrule
1 & $\mathrm{Uniform}\sim\xi\in\Xi=[a,b]$ & $\mathbf{E}[\xi]=400$ N/m & Gauss-Legendre (100 points)\\
& $a=340$ N/m, $b=460$ N/m & $\mathrm{Var}[\xi]=1\,200$ N$^2$/m$^2$ & \\
\midrule
2 & $\mathrm{Beta}(\alpha,\beta)\sim\xi\in\Xi=[a,b]$ & $\mathbf{E}[\xi]\approx374.3$ N/m & Gauss-Jacobi (95 points)\\
& $\alpha=2$, $\beta=5$ & $\mathrm{Var}[\xi]\approx367.3$ N$^2$/m$^2$ & $\alpha_J:=\beta-1=4$ \\
& $a=340$ N/m, $b=460$ N/m & & $\beta_J:=\alpha-1=1$ \\
\midrule
3 & $\mathrm{Gamma}(\alpha,\beta)\sim\xi\in\Xi=[a,\infty)$ & $\mathbf{E}[\xi]=440$ N/m & Gauss-Laguerre (135 points)\\
& $\alpha=10$, $\beta=1/10$ & $\mathrm{Var}[\xi]=1\,000$ N$^2$/m$^2$ & $\alpha_L:=\alpha-1=9$ \\
& $a=340$ N/m & & \\
\midrule
\multicolumn{4}{@{}p{\textwidth}@{}}{\footnotesize $^*$Probability density function, $f:\Xi\to\mathbb{R}^+_0$:%
\begin{equation*}\mathrm{Uniform}\sim f(\xi)=\frac{1}{b-a},\quad\mathrm{Beta}\sim f(\xi)=\frac{(\xi-a)^{\alpha-1}\,(b-\xi)^{\beta-1}}{(b-a)^{\alpha+\beta-1}\,\mathrm{B}(\alpha,\beta)}\quad\text{and}\quad\mathrm{Gamma}\sim f(\xi)=\frac{\beta^\alpha}{\Gamma(\alpha)}(\xi-a)^{\alpha-1}\exp(-\beta\,(\xi-a)).%
\end{equation*}}\\[-1ex]
\bottomrule
\end{tabular}
\label{tab1CaseStudies10}
\end{table}

\begin{figure}
\centering
\begin{subfigure}[b]{0.495\textwidth}
\includegraphics[width=\textwidth]{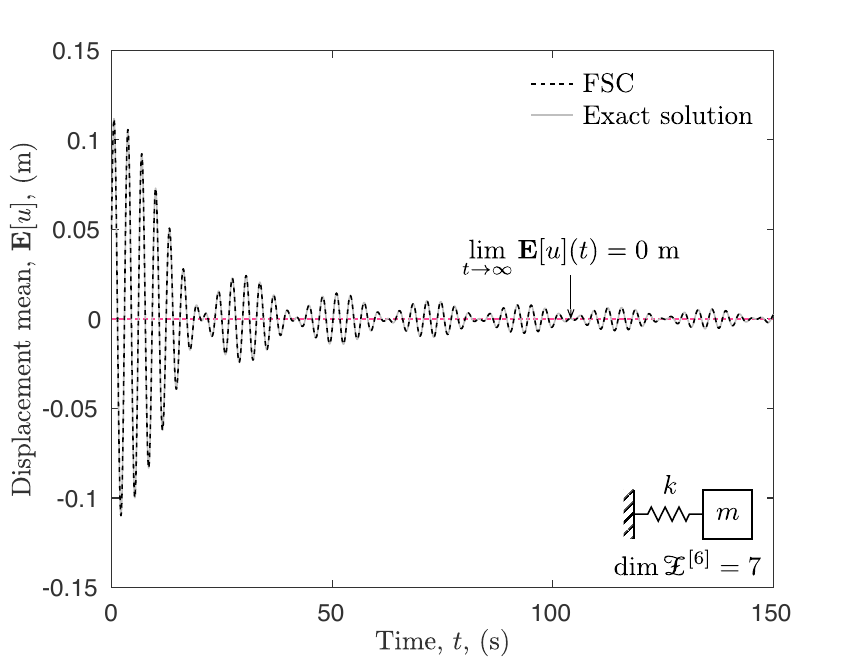}
\caption{Mean}
\label{fig1SDOF1_Uniform_FSC_Disp_Mean_7}
\end{subfigure}\hfill
\begin{subfigure}[b]{0.495\textwidth}
\includegraphics[width=\textwidth]{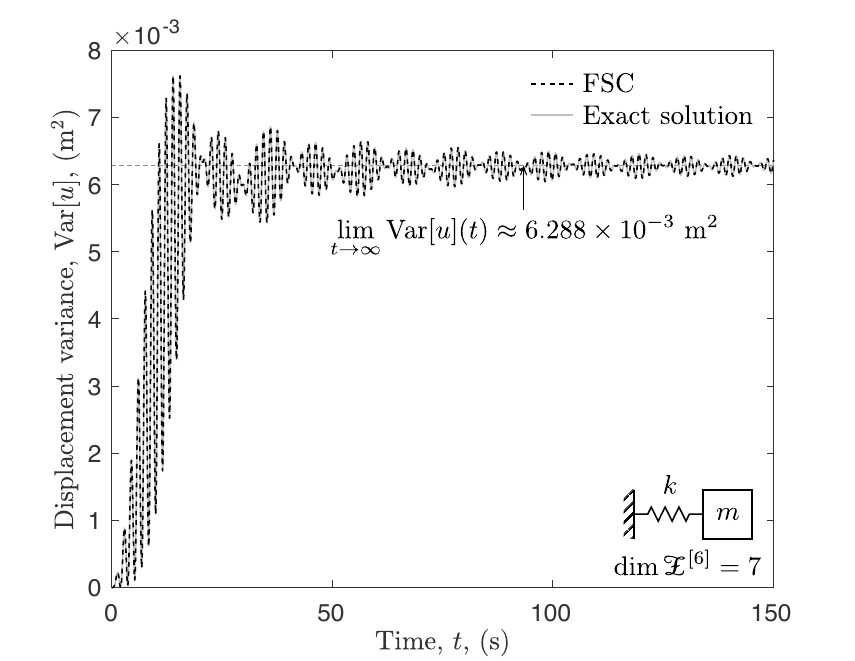}
\caption{Variance}
\label{fig1SDOF1_Uniform_FSC_Disp_Var_7}
\end{subfigure}
\caption{Evolution of $\mathbf{E}[u]$ and $\mathrm{Var}[u]$ for the case when the $p$-discretization level of RFS is $\mathscr{Z}^{[6]}$ and $\mu\sim\mathrm{Uniform}$}
\label{fig1SDOF1_Uniform_FSC_Disp_7}
\end{figure}

\begin{figure}
\centering
\begin{subfigure}[b]{0.495\textwidth}
\includegraphics[width=\textwidth]{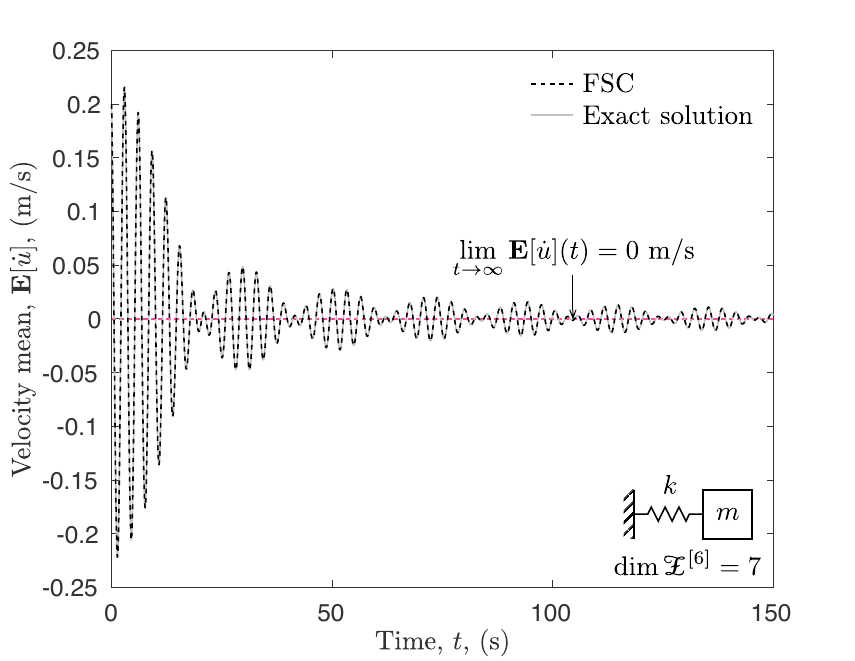}
\caption{Mean}
\label{fig1SDOF1_Uniform_FSC_Vel_Mean_7}
\end{subfigure}\hfill
\begin{subfigure}[b]{0.495\textwidth}
\includegraphics[width=\textwidth]{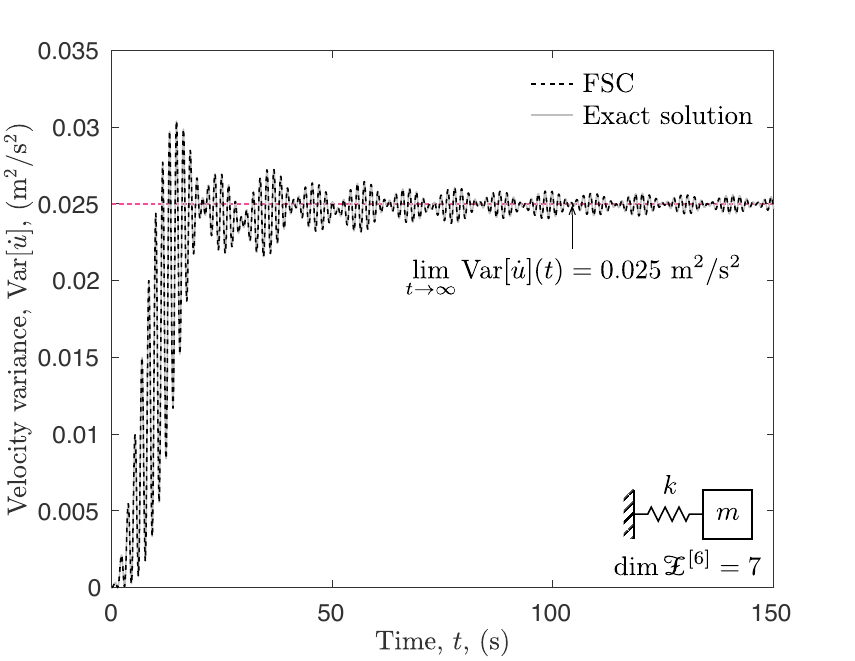}
\caption{Variance}
\label{fig1SDOF1_Uniform_FSC_Vel_Var_7}
\end{subfigure}
\caption{Evolution of $\mathbf{E}[\dot{u}]$ and $\mathrm{Var}[\dot{u}]$ for the case when the $p$-discretization level of RFS is $\mathscr{Z}^{[6]}$ and $\mu\sim\mathrm{Uniform}$}
\label{fig1SDOF1_Uniform_FSC_Vel_7}
\end{figure}

Figs.~\ref{fig1SDOF1_Uniform_FSC_Disp_7} and \ref{fig1SDOF1_Uniform_FSC_Vel_7} show the evolution of the mean and variance of the system's state.
From these figures, it is observed that the response obtained with FSC using only 7 basis vectors has the ability to reproduce the exact response (from Appendix \ref{appsec1SDOF}) with high fidelity.
This is the reason why the two plots appear to be indistinguishable from each other.
The figures also show the limit values for each response computed using the exact expressions given by \eqref{appeq1SDOF5500} and \eqref{appeq1SDOF7000}.

\begin{figure}
\centering
\begin{subfigure}[b]{0.495\textwidth}
\includegraphics[width=\textwidth]{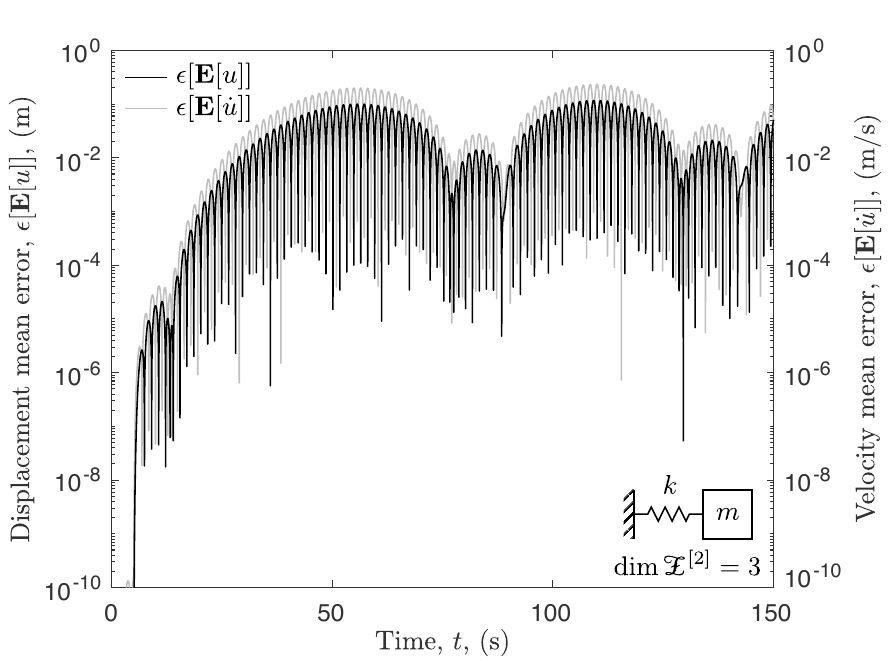}
\caption{Mean error for $\mathscr{Z}^{[2]}$}
\label{fig1SDOF1_Uniform_FSC_Mean_3_Error}
\end{subfigure}\hfill
\begin{subfigure}[b]{0.495\textwidth}
\includegraphics[width=\textwidth]{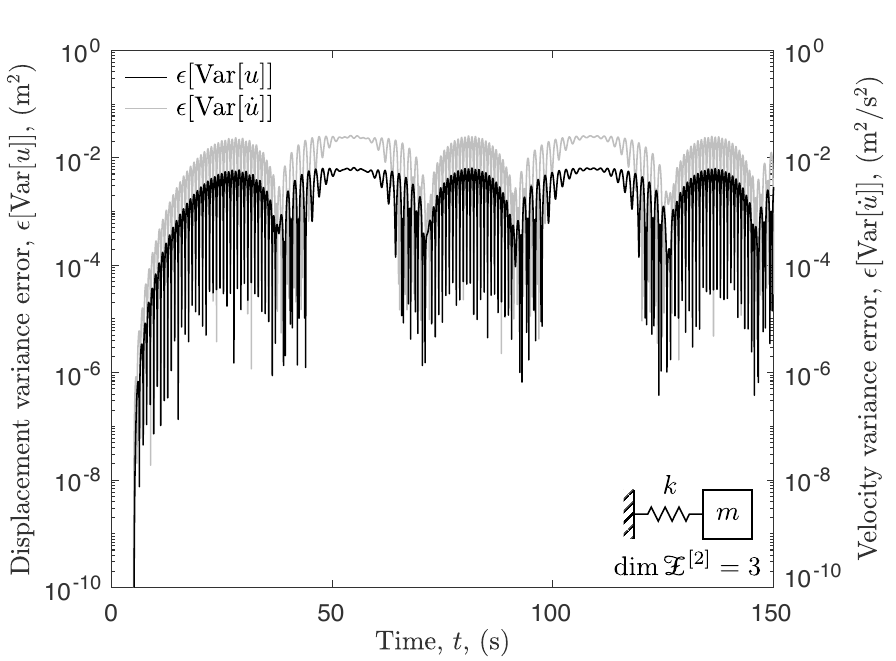}
\caption{Variance error for $\mathscr{Z}^{[2]}$}
\label{fig1SDOF1_Uniform_FSC_Var_3_Error}
\end{subfigure}\quad
\begin{subfigure}[b]{0.495\textwidth}
\includegraphics[width=\textwidth]{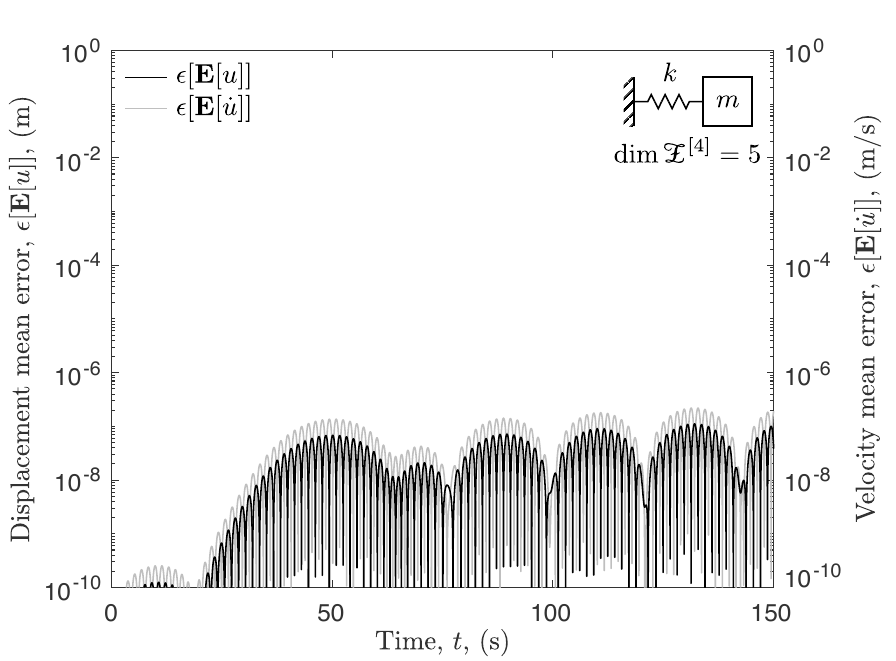}
\caption{Mean error for $\mathscr{Z}^{[4]}$}
\label{fig1SDOF1_Uniform_FSC_Mean_5_Error}
\end{subfigure}\hfill
\begin{subfigure}[b]{0.495\textwidth}
\includegraphics[width=\textwidth]{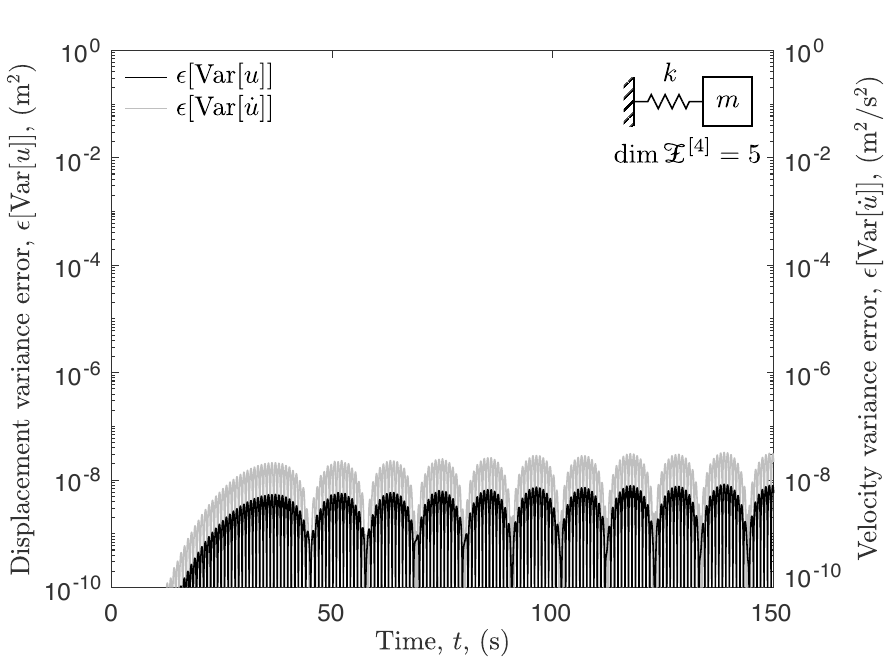}
\caption{Variance error for $\mathscr{Z}^{[4]}$}
\label{fig1SDOF1_Uniform_FSC_Var_5_Error}
\end{subfigure}\quad
\begin{subfigure}[b]{0.495\textwidth}
\includegraphics[width=\textwidth]{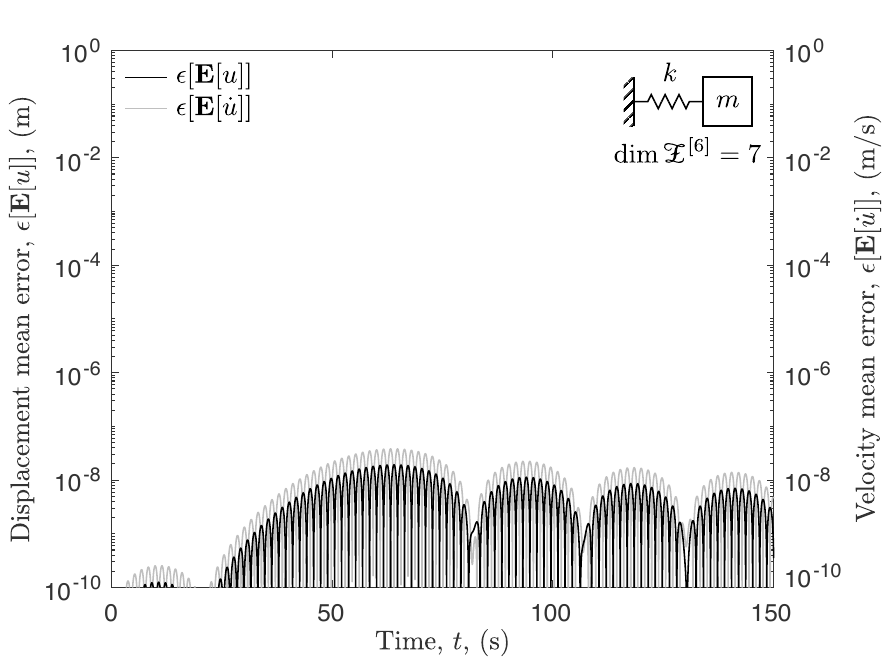}
\caption{Mean error for $\mathscr{Z}^{[6]}$}
\label{fig1SDOF1_Uniform_FSC_Mean_7_Error}
\end{subfigure}\hfill
\begin{subfigure}[b]{0.495\textwidth}
\includegraphics[width=\textwidth]{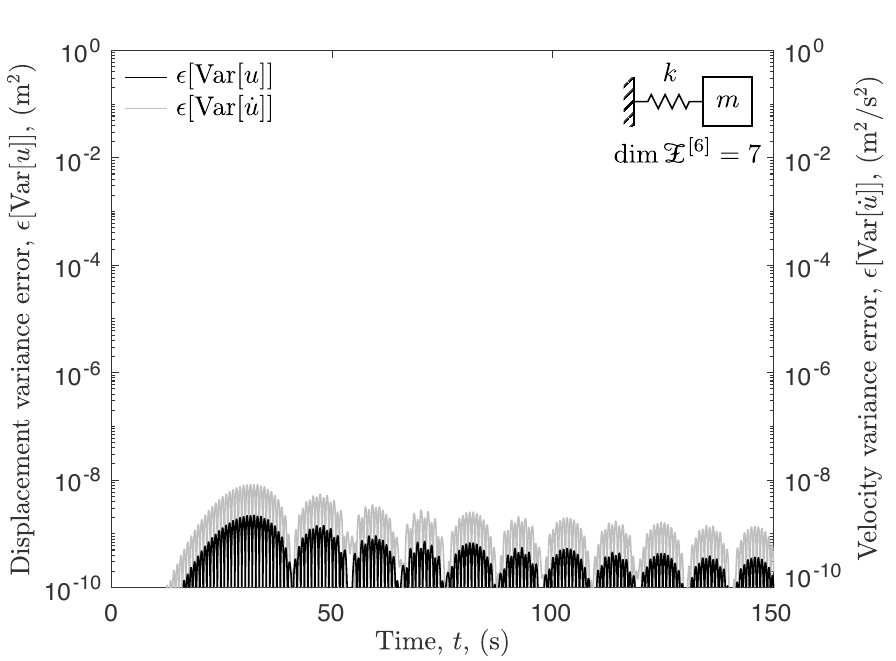}
\caption{Variance error for $\mathscr{Z}^{[6]}$}
\label{fig1SDOF1_Uniform_FSC_Var_7_Error}
\end{subfigure}
\caption{Local error evolution of $\mathbf{E}[u]$, $\mathrm{Var}[u]$, $\mathbf{E}[\dot{u}]$ and $\mathrm{Var}[\dot{u}]$ for different $p$-discretization levels of RFS and for $\mu\sim\mathrm{Uniform}$}
\label{fig1SDOF1_Uniform_FSC_Error}
\end{figure}

\begin{figure}
\centering
\begin{subfigure}[b]{0.495\textwidth}
\includegraphics[width=\textwidth]{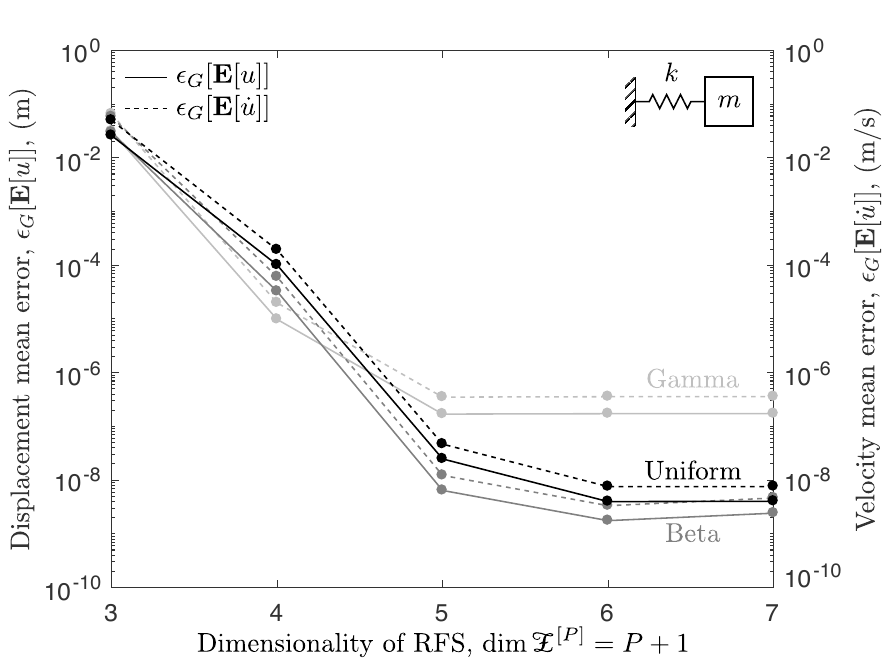}
\caption{Mean error}
\label{fig1SDOF1_FSC_Mean_GlobalError}
\end{subfigure}\hfill
\begin{subfigure}[b]{0.495\textwidth}
\includegraphics[width=\textwidth]{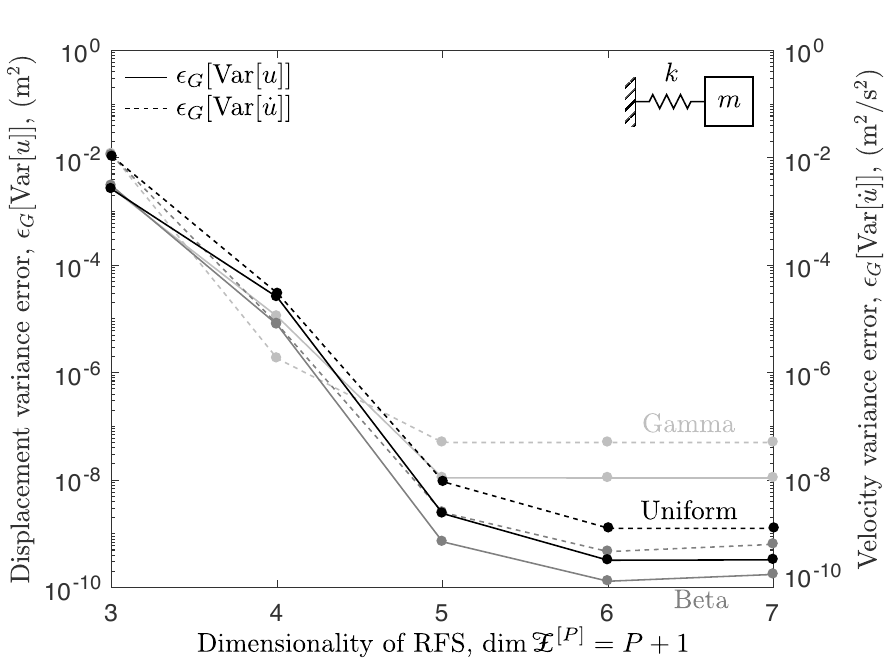}
\caption{Variance error}
\label{fig1SDOF1_FSC_Var_GlobalError}
\end{subfigure}
\caption{Global error of $\mathbf{E}[u]$, $\mathrm{Var}[u]$, $\mathbf{E}[\dot{u}]$ and $\mathrm{Var}[\dot{u}]$ for different $p$-discretization levels of RFS}
\label{fig1SDOF1_FSC_GlobalError}
\end{figure}

\begin{figure}
\centering
\begin{subfigure}[b]{0.495\textwidth}
\includegraphics[width=\textwidth]{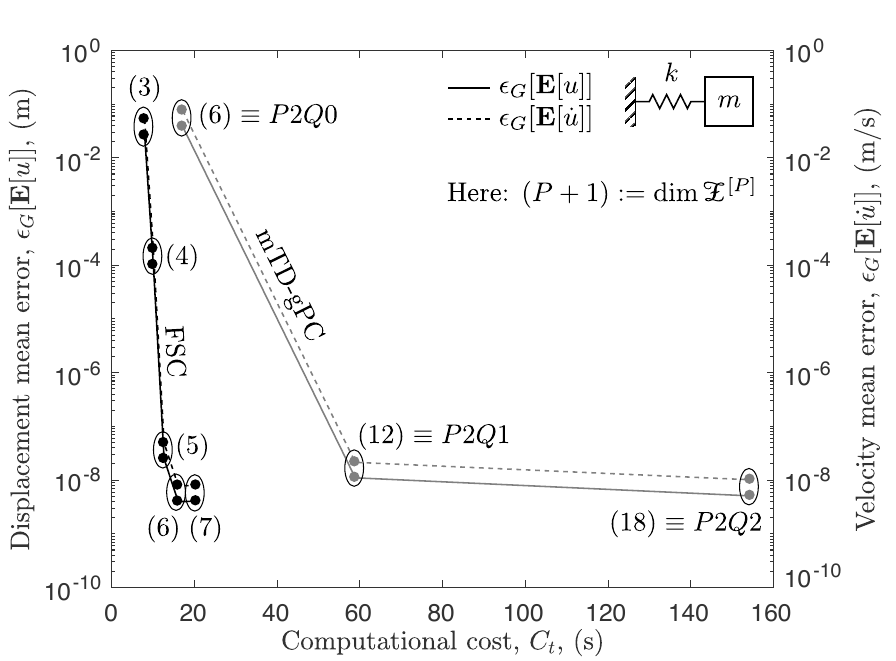}
\caption{Mean error}
\label{fig1SDOF1_Uniform_ComputationalCost_Mean}
\end{subfigure}\hfill
\begin{subfigure}[b]{0.495\textwidth}
\includegraphics[width=\textwidth]{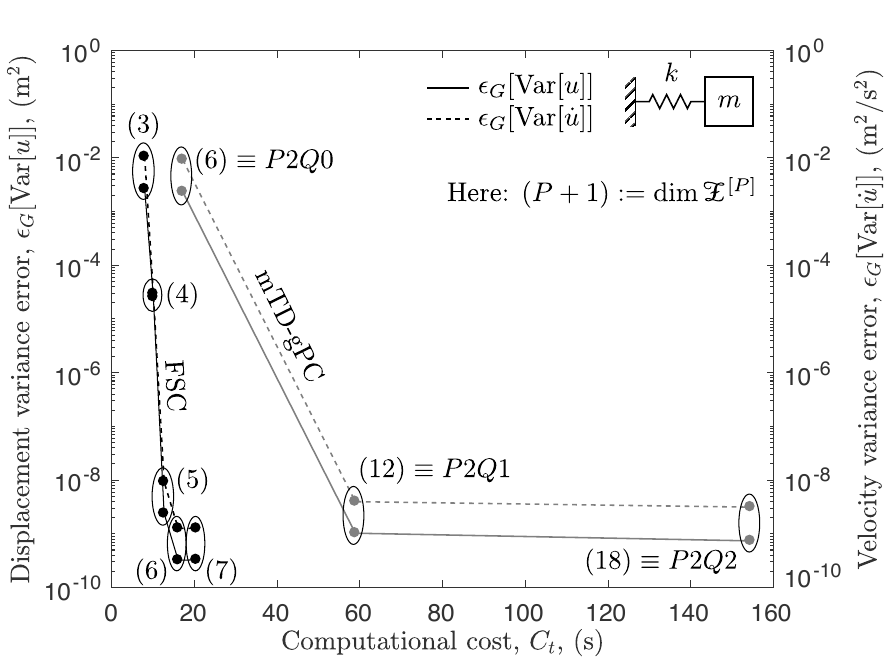}
\caption{Variance error}
\label{fig1SDOF1_Uniform_ComputationalCost_Var}
\end{subfigure}
\caption{Global error versus computational cost for $\mu\sim\mathrm{Uniform}$}
\label{fig1SDOF1_Uniform_ComputationalCost}
\end{figure}

Fig.~\ref{fig1SDOF1_Uniform_FSC_Error} presents the local errors in mean and variance of the system's state using different choices of number of basis vectors ranging from 3 to 7.
Note that even though the FSC scheme requires $P$ to be greater than or equal to 3, here we also study the case when $P=2$ for sake of comparison.
Fig.~\ref{fig1SDOF1_Uniform_FSC_Error} shows that as the number of basis vectors increases, so does the accuracy of the results.
In particular, by increasing the number of basis vectors from 3 to 5, the accuracy of the results improves significantly by an order of magnitude of 6 (approximately from $10^{-1}$ to $10^{-7}$ for the mean), whereas when the number of basis vectors is increased from 5 to 7, the improvement in error is more moderate (to approximately $10^{-8}$ for the mean).
However, for the 7-basis-vector case we do see an improvement in the accuracy of the solution as time progresses due to the increase in the number of basis vectors used.
Fig.~\ref{fig1SDOF1_FSC_GlobalError} presents the convergence of global errors as a function of the number of basis vectors and the different distributions used to define the stochasticity of $k$.
The FSC scheme achieves exponential convergence when 3, 4 and 5 basis vectors are used, but adding more basis vectors does not improve the accuracy of the response.
The primary reason for this slow-down in convergence is that the accuracy of the solution is limited by machine precision and the fact that the probability information is being transferred in the mean-square sense at every time step (for $P$ is finite).
In fact, these plots indicate that there is no reason to implement more than 5 basis vectors into the simulation as it does not improve the accuracy of the results significantly.
It is also apparent from this figure that when $k$ is assumed gamma-distributed, the results are not as accurate as those obtained from the uniform and beta distributions.
The reason behind this outcome is that when $k$ is gamma-distributed its support is unbounded, which from a numerical viewpoint leads to the dreaded case of unbounded basis vectors.

Fig.~\ref{fig1SDOF1_Uniform_ComputationalCost} plots the global errors as a function of computational cost\footnote{All problems in this work were run in MATLAB R2016b \cite{matlabr2016b} on a 2017 MacBook Pro with quad-core \emph{3.1 GHz Intel Core i7} processor (hyper-threading technology enabled), \emph{16 GB 2133 MHz LPDDR3} memory and \emph{1 TB PCI-Express SSD} storage (APFS-formatted), running macOS Mojave (version 10.14.6).} of the FSC and mTD-gPC schemes expressed in terms of the wall-clock time taken to complete the computation.
The implementation of both schemes was optimized as much as reasonably possible, and the labels $P2Q0$, $P2Q1$ and $P2Q2$ are defined in Ref.~\cite{heuveline2014hybrid} (Pg.~45).
The comparison of computational cost is shown here only for Case 1 of Table \ref{tab1CaseStudies10} (for sake of brevity), but similar trends are observed for cases 2 and 3 as well.
Note that FSC is much faster in comparison to mTD-gPC for a similar level of error.
For instance, in order to attain a global error of approximately $10^{-8}$, FSC runs about 3.5 times faster than mTD-gPC.
This is because, in general, FSC requires much fewer basis vectors than mTD-gPC---to achieve an error of about $10^{-8}$, FSC requires only 6 basis vectors in comparison to 12 for mTD-gPC.
Another reason for the superior efficiency of FSC is that, for mTD-gPC, the orthogonalization process needs to be conducted three times when the random basis is demanded to be updated during the simulation (one time for the monomials of $u$, another time for the monomials of $\dot{u}$, and one more time after performing the tensor product)\footnote{We tested the method by orthogonalizing only once---namely, after performing the tensor product between the monomials of $u$ and $\dot{u}$---, and found that the accuracy of the results degrades noticeably.}.
This also explains why using 6 basis vectors in both methods, FSC runs slightly faster than mTD-gPC.
Furthermore, we see that the probability information is better encoded in FSC because it uses a fewer number of basis vectors to achieve the same level of accuracy.
Fig.~\ref{fig1SDOF1_Uniform_ComputationalCost} also reveals that increasing the number of basis vectors from 6 to 7 for FSC and 12 to 18 for mTD-gPC does not improve the accuracy of the results significantly.
This, again, is because of the limited precision of the machine and the fact that the probability information is being transferred in the mean-square sense.

\subsection{Single-degree-of-freedom system under forced vibration}\label{sec1NumRes20}

In this example, we show that the number of basis vectors needed in the simulation does not increase when the dimensionality of the random space increases.
For this, we consider the same system described in the Example \ref{sec1NumRes10} (including the same deterministic initial conditions), with the only difference being that the system is subjected to a stochastic external force given by $p(t,\cdot\,)=q\sin(t)$.
That is, the system is now governed by $m\ddot{u}+ku=p$.
Here the stiffness $k(\xi)=\xi^1$ is taken to be the same as Case 1 of Table \ref{tab1CaseStudies10} with $\bar{\Xi}_1=[340,460]$ N/m.
The amplitude of the external force $q(\xi)=\xi^2$ is assumed beta-distributed with parameters $\alpha=2$ and $\beta=5$ in $\bar{\Xi}_2=[51,69]$ N, giving thereby probability moments: $\mathbf{E}[\xi^2]\approx56.14$ N and $\mathrm{Var}[\xi^2]\approx8.265$ N$^2$.
Because two random variables are present in the mathematical model, the random domain of the system is 2-dimensional, and thus, it is defined by $\Xi=\bar{\Xi}_1\times\bar{\Xi}_2$ with $\mu\sim\mathrm{Uniform}\otimes\mathrm{Beta}$.
For this example, the inner products are computed using a quadrature rule constructed by performing a cartesian product between 100 Gauss-Legendre points distributed along the $\xi^1$-axis and 95 Gauss-Jacobi points distributed along the $\xi^2$-axis.
Finally, the gPC method (with $P=8$) is used for the first 0.5 seconds of the simulation to allow the stochasticity of the system to develop sufficiently before using the FSC scheme.

\begin{figure}
\centering
\begin{subfigure}[b]{0.495\textwidth}
\includegraphics[width=\textwidth]{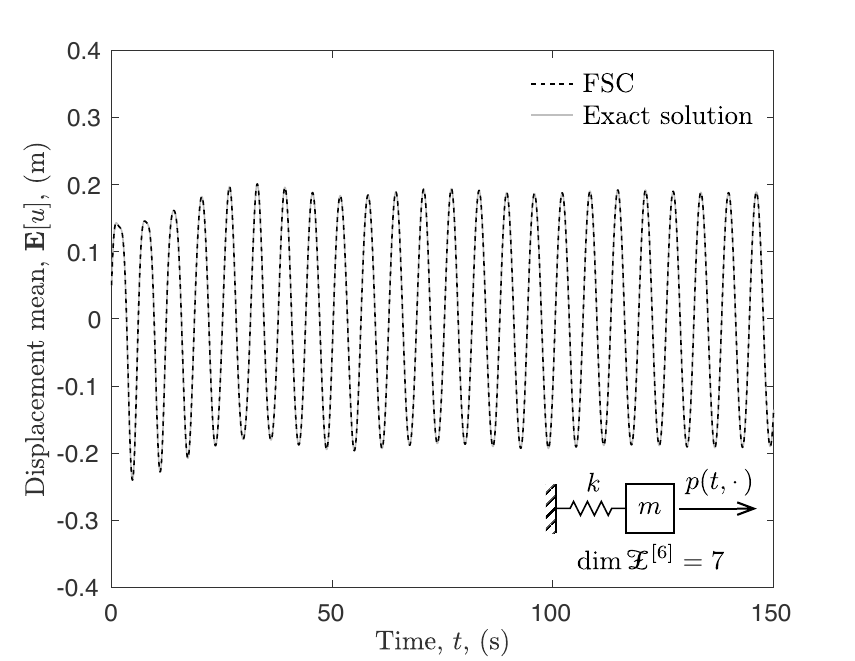}
\caption{Mean}
\label{fig1SDOF2_UniformBeta_FSC_Disp_Mean_7}
\end{subfigure}\hfill
\begin{subfigure}[b]{0.495\textwidth}
\includegraphics[width=\textwidth]{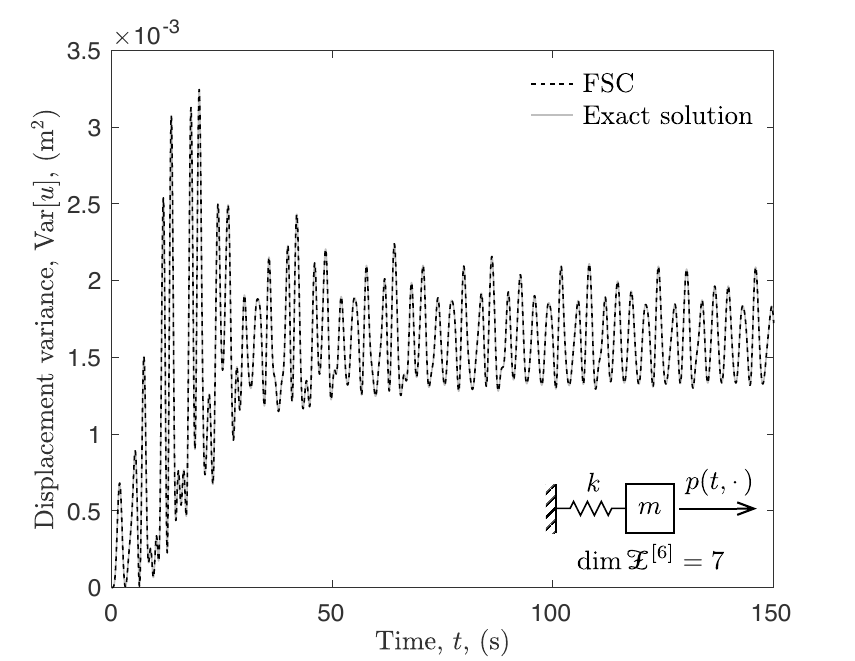}
\caption{Variance}
\label{fig1SDOF2_UniformBeta_FSC_Disp_Var_7}
\end{subfigure}
\caption{Evolution of $\mathbf{E}[u]$ and $\mathrm{Var}[u]$ for the case when the $p$-discretization level of RFS is $\mathscr{Z}^{[6]}$ and $\mu\sim\mathrm{Uniform}\otimes\mathrm{Beta}$}
\label{fig1SDOF2_UniformBeta_FSC_Disp_7}
\end{figure}

\begin{figure}
\centering
\begin{subfigure}[b]{0.495\textwidth}
\includegraphics[width=\textwidth]{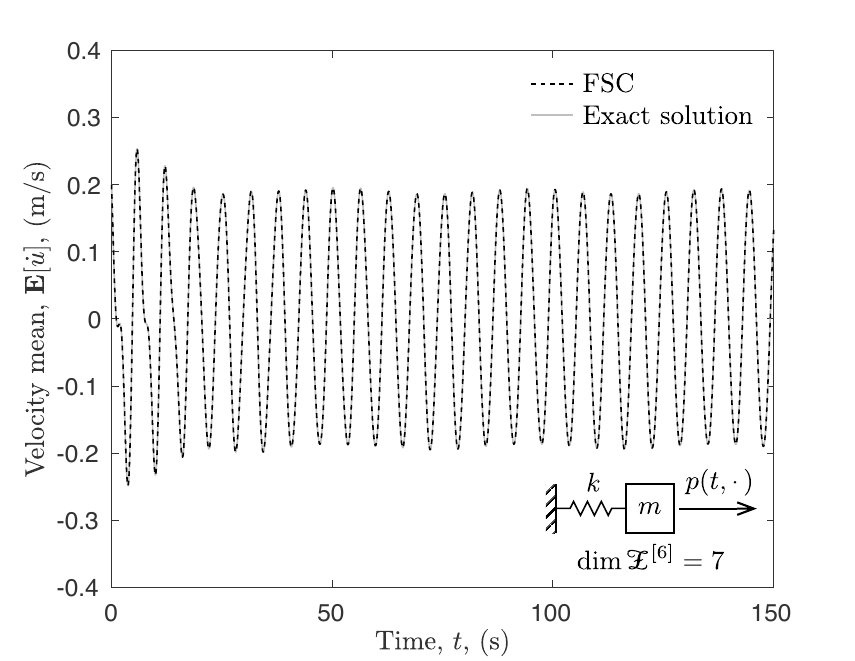}
\caption{Mean}
\label{fig1SDOF2_UniformBeta_FSC_Vel_Mean_7}
\end{subfigure}\hfill
\begin{subfigure}[b]{0.495\textwidth}
\includegraphics[width=\textwidth]{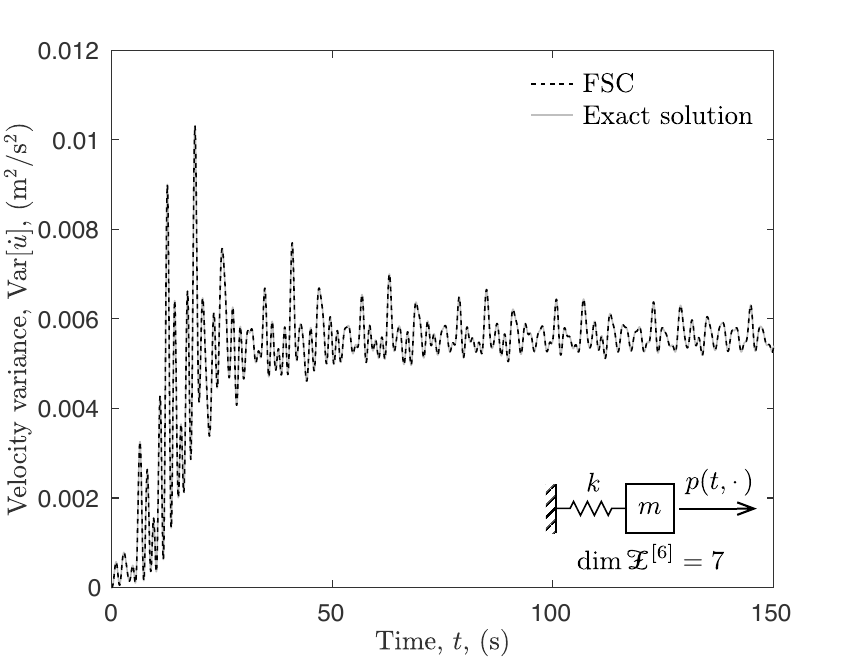}
\caption{Variance}
\label{fig1SDOF2_UniformBeta_FSC_Vel_Var_7}
\end{subfigure}
\caption{Evolution of $\mathbf{E}[\dot{u}]$ and $\mathrm{Var}[\dot{u}]$ for the case when the $p$-discretization level of RFS is $\mathscr{Z}^{[6]}$ and $\mu\sim\mathrm{Uniform}\otimes\mathrm{Beta}$}
\label{fig1SDOF2_UniformBeta_FSC_Vel_7}
\end{figure}

\begin{figure}
\centering
\begin{subfigure}[b]{0.495\textwidth}
\includegraphics[width=\textwidth]{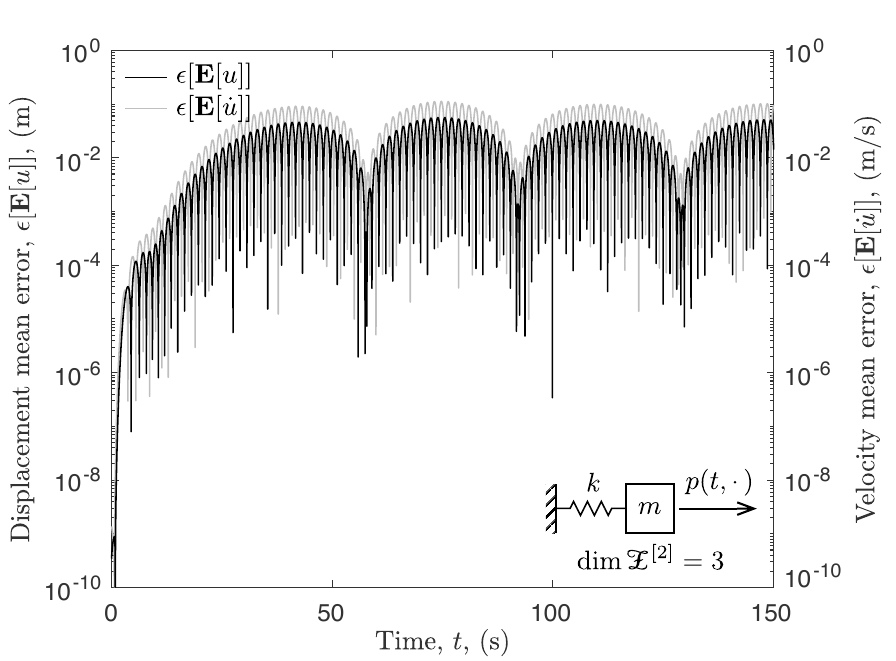}
\caption{Mean error for $\mathscr{Z}^{[2]}$}
\label{fig1SDOF2_UniformBeta_FSC_Mean_3_Error}
\end{subfigure}\hfill
\begin{subfigure}[b]{0.495\textwidth}
\includegraphics[width=\textwidth]{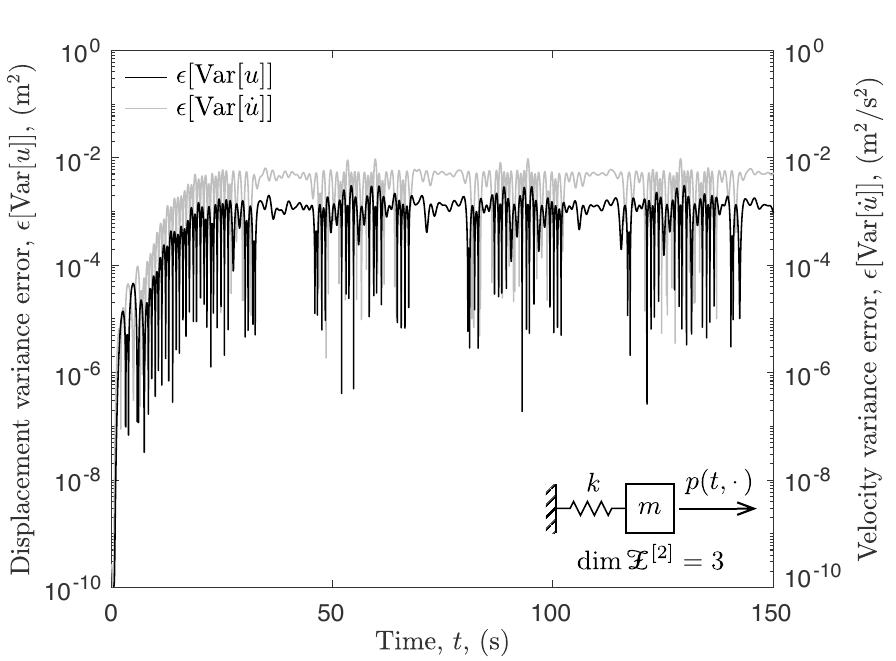}
\caption{Variance error for $\mathscr{Z}^{[2]}$}
\label{fig1SDOF2_UniformBeta_FSC_Var_3_Error}
\end{subfigure}\quad
\begin{subfigure}[b]{0.495\textwidth}
\includegraphics[width=\textwidth]{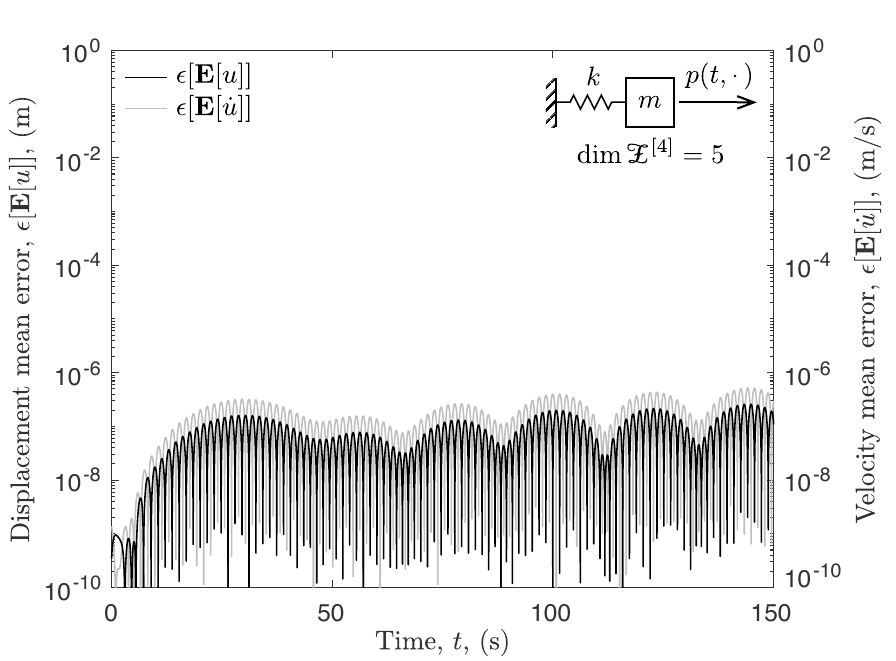}
\caption{Mean error for $\mathscr{Z}^{[4]}$}
\label{fig1SDOF2_UniformBeta_FSC_Mean_5_Error}
\end{subfigure}\hfill
\begin{subfigure}[b]{0.495\textwidth}
\includegraphics[width=\textwidth]{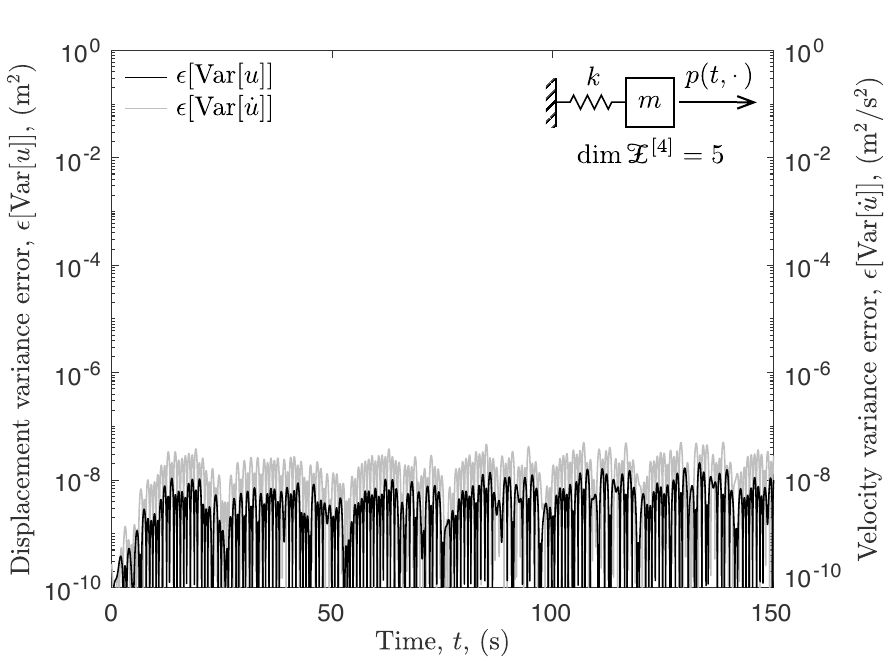}
\caption{Variance error for $\mathscr{Z}^{[4]}$}
\label{fig1SDOF2_UniformBeta_FSC_Var_5_Error}
\end{subfigure}\quad
\begin{subfigure}[b]{0.495\textwidth}
\includegraphics[width=\textwidth]{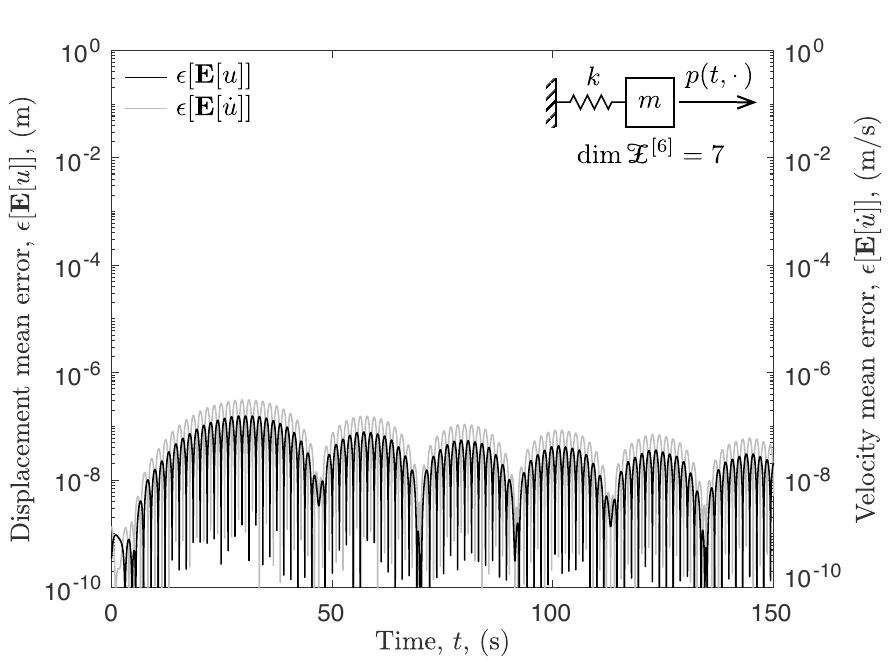}
\caption{Mean error for $\mathscr{Z}^{[6]}$}
\label{fig1SDOF2_UniformBeta_FSC_Mean_7_Error}
\end{subfigure}\hfill
\begin{subfigure}[b]{0.495\textwidth}
\includegraphics[width=\textwidth]{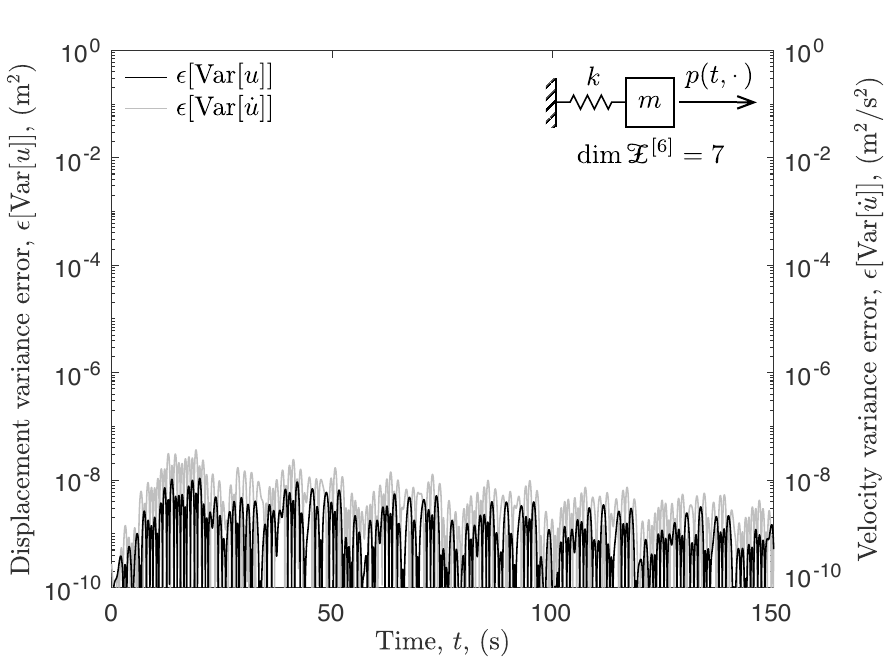}
\caption{Variance error for $\mathscr{Z}^{[6]}$}
\label{fig1SDOF2_UniformBeta_FSC_Var_7_Error}
\end{subfigure}
\caption{Local error evolution of $\mathbf{E}[u]$, $\mathrm{Var}[u]$, $\mathbf{E}[\dot{u}]$ and $\mathrm{Var}[\dot{u}]$ for different $p$-discretization levels of RFS and for $\mu\sim\mathrm{Uniform}\otimes\mathrm{Beta}$}
\label{fig1SDOF2_UniformBeta_FSC_Error}
\end{figure}

\begin{figure}
\centering
\begin{subfigure}[b]{0.495\textwidth}
\includegraphics[width=\textwidth]{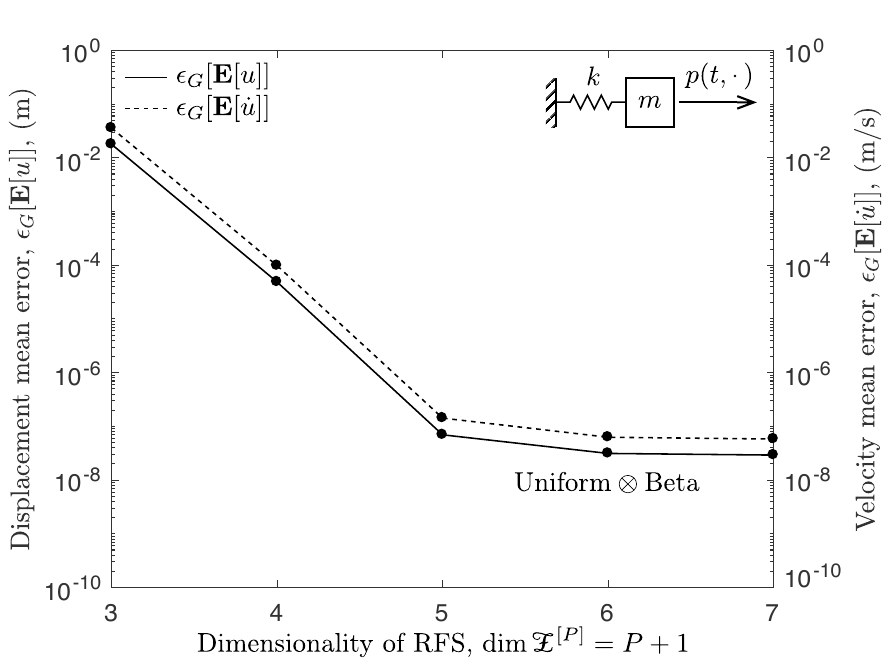}
\caption{Mean error}
\label{fig1SDOF2_FSC_Mean_GlobalError}
\end{subfigure}\hfill
\begin{subfigure}[b]{0.495\textwidth}
\includegraphics[width=\textwidth]{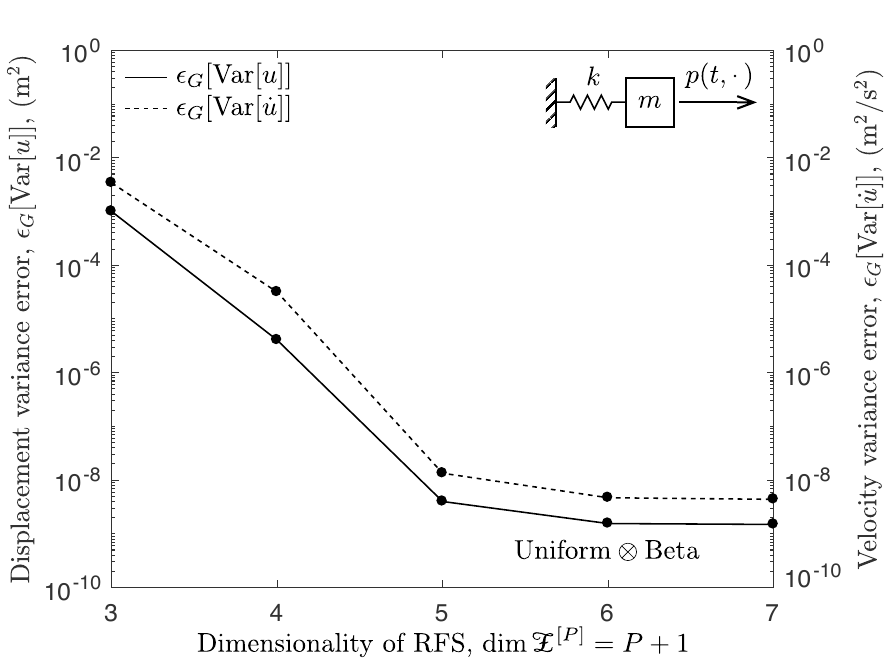}
\caption{Variance error}
\label{fig1SDOF2_FSC_Var_GlobalError}
\end{subfigure}
\caption{Global error of $\mathbf{E}[u]$, $\mathrm{Var}[u]$, $\mathbf{E}[\dot{u}]$ and $\mathrm{Var}[\dot{u}]$ for different $p$-discretization levels of RFS}
\label{fig1SDOF2_FSC_GlobalError}
\end{figure}

\begin{figure}
\centering
\begin{subfigure}[b]{0.495\textwidth}
\includegraphics[width=\textwidth]{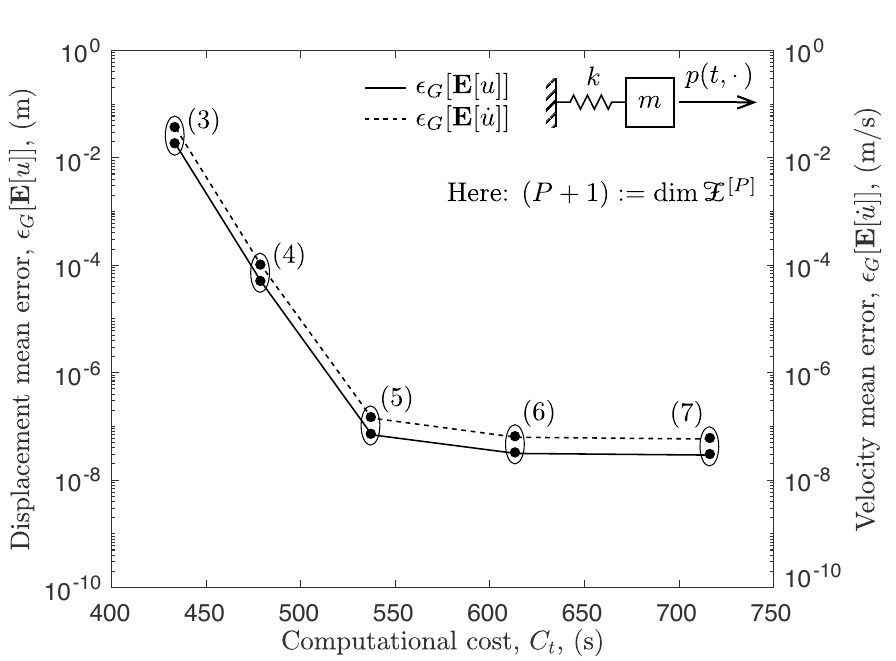}
\caption{Mean error}
\label{fig1SDOF2_UniformBeta_ComputationalCost_Mean}
\end{subfigure}\hfill
\begin{subfigure}[b]{0.495\textwidth}
\includegraphics[width=\textwidth]{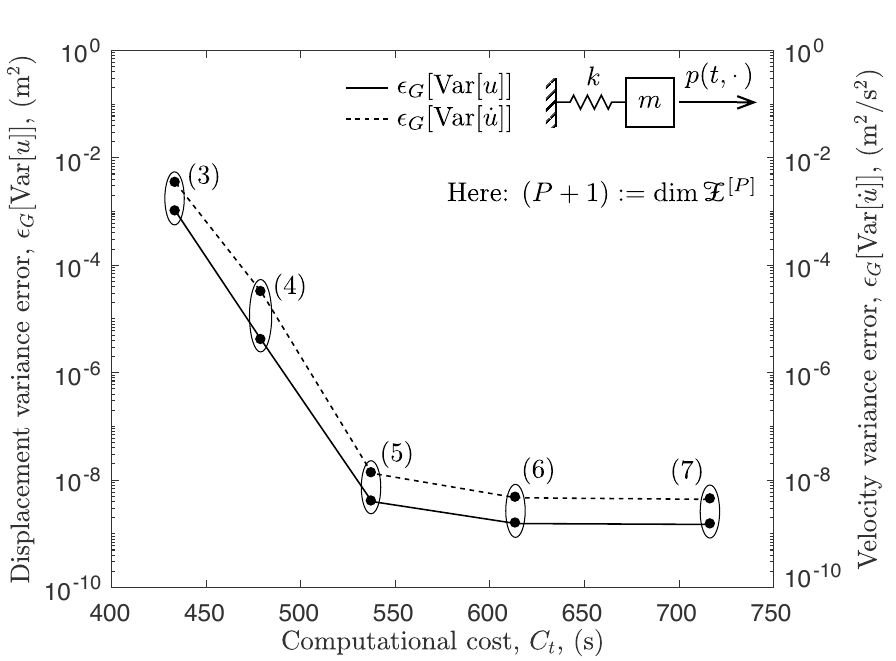}
\caption{Variance error}
\label{fig1SDOF2_UniformBeta_ComputationalCost_Var}
\end{subfigure}
\caption{Global error versus computational cost for $\mu\sim\mathrm{Uniform}\otimes\mathrm{Beta}$}
\label{fig1SDOF2_UniformBeta_ComputationalCost}
\end{figure}

Figs.~\ref{fig1SDOF2_UniformBeta_FSC_Disp_7} and \ref{fig1SDOF2_UniformBeta_FSC_Vel_7} show the evolution of the mean and variance of the system's state.
As in the previous example, the numerical solution obtained using FSC with 7 basis vectors is indistinguishable from the exact response\footnote{To obtain the `exact' solution for $\mathbf{E}[u]$, $\mathrm{Var}[u]$, $\mathbf{E}[\dot{u}]$ and $\mathrm{Var}[\dot{u}]$, the corresponding values for $\mathbf{E}[u(t_i,\cdot\,)]$, $\mathrm{Var}[u(t_i,\cdot\,)]$, $\mathbf{E}[\dot{u}(t_i,\cdot\,)]$ and $\mathrm{Var}[\dot{u}(t_i,\cdot\,)]$ were computed at each instant of time $t_i\in\mathfrak{T}$ using the {\tt vpaintegral} command (provided in the MATLAB's Symbolic Math Toolbox \cite{matlabr2016b}) with {\tt RelTol} set equal to $10^{-14}$. The exact displacement response, $u$, is well known and can be found in any structural dynamics textbook, e.g.~\cite{chopra2012dynamics,humar2012dynamics}.}.
Figs.~\ref{fig1SDOF2_UniformBeta_FSC_Error} and \ref{fig1SDOF2_FSC_GlobalError} depict the local and global errors in mean and variance of the system's state.
Here we also notice the same trend found in Figs.~\ref{fig1SDOF1_Uniform_FSC_Error} and \ref{fig1SDOF1_FSC_GlobalError}.
That is, as the number of basis vectors increases, so does the accuracy of the results.
Moreover, when the number of basis vectors increases from 3 to 5, the error for the mean drops down from approximately $10^{-1}$ to $10^{-7}$, whereas increasing the number of basis vectors from 5 to 7 does not result in a noticeable improvement in the computation of the probability moments.
Note that using the same number of basis vectors as in Example \ref{sec1NumRes10} led to similar levels of error, even though the dimensionality of the random space in this example is twice that of Example \ref{sec1NumRes10}.
However, despite the number of basis vectors being the same, the computational cost of this example is much higher than that of Example \ref{sec1NumRes10} (as per Fig.~\ref{fig1SDOF2_UniformBeta_ComputationalCost}) because of the increase in the number of quadrature points needed to compute the inner products accurately.

\subsection{Nonlinear single-degree-of-freedom system under free vibration}\label{sec1NumRes30}

In this last example, we explore the nonlinear behavior of a single-degree-of-freedom system in order to test the ability of the FSC method to solve nonlinear problems.
The governing differential equation for this system is set to be given by
\begin{equation}\label{eq1NumRes300}
m\ddot{u}+(1+\rho u^2)ku=0,
\end{equation}
where $m=100$ kg is the mass of the system, $k(\xi)=\xi^1$ is a beta-distributed random variable representing the strength of the stiffness which is given by Case 2 of Table \ref{tab1CaseStudies10} with $\bar{\Xi}_1=[340,460]$ N/m, and $\rho(\xi)=\xi^2$ is a uniformly-distributed random variable denoting the contributing factor to the nonlinearity of the system in $\bar{\Xi}_2=[-20,-30]$ m$^{-2}$.
The probability moments for $\xi^2$ are thus: $\mathbf{E}[\xi^2]=-25$ m$^{-2}$ and $\mathrm{Var}[\xi^2]\approx8.333$ m$^{-4}$.
The system has an initial displacement of $u(0,\cdot\,)\equiv0.05$ m and an initial velocity of $\dot{u}(0,\cdot\,)\equiv0.20$ m/s.
As in the previous example, the random domain of the system is two-dimensional and defined by $\Xi=\bar{\Xi}_1\times\bar{\Xi}_2$ with $\mu\sim\mathrm{Beta}\otimes\mathrm{Uniform}$.
The inner products are again computed with a Gaussian quadrature rule using the same number of points indicated in the previous section, and the gPC method (with $P=8$) is used for the first second of the simulation.
The system is integrated over time using the RK4 method with a time-step size of $\Delta t=0.005$ s, and the simulation is set to last $T=150$ s.

\begin{remark}
According to \eqref{eq1SolSpeApp1060}, the temporal function $\mathcal{F}^i[u^j,\dot{u}^j]$ associated with $\mathcal{F}[u,\dot{u}]=(1+\rho u^2)ku$ is given by:
\begin{equation*}
\mathcal{F}^i[u^j,\dot{u}^j](t)=\frac{\langle\Psi_i,k\Psi_j\rangle}{\langle\Psi_i,\Psi_i\rangle}\,u^j(t)
+\frac{\langle\Psi_i,\rho k\Psi_j\Psi_k\Psi_l\rangle}{\langle\Psi_i,\Psi_i\rangle}\,u^j(t)\,u^k(t)\,u^l(t).
\end{equation*}
Note that a summation sign is implied over every repeated index $j$, $k$ and $l$.
\end{remark}

\begin{figure}
\centering
\begin{subfigure}[b]{0.495\textwidth}
\includegraphics[width=\textwidth]{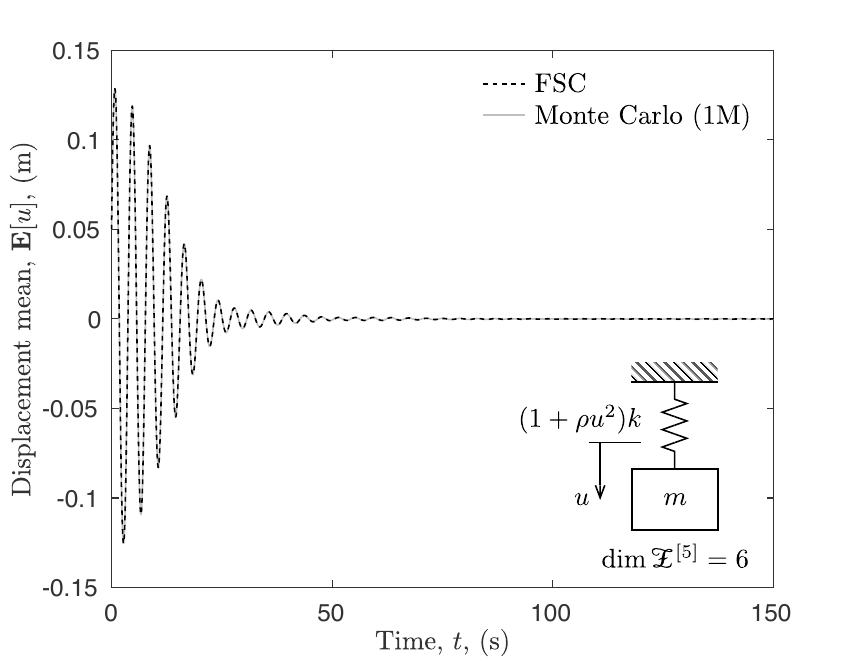}
\caption{Mean}
\label{fig1SDOF3_BetaUniform_FSC_Disp_Mean_6}
\end{subfigure}\hfill
\begin{subfigure}[b]{0.495\textwidth}
\includegraphics[width=\textwidth]{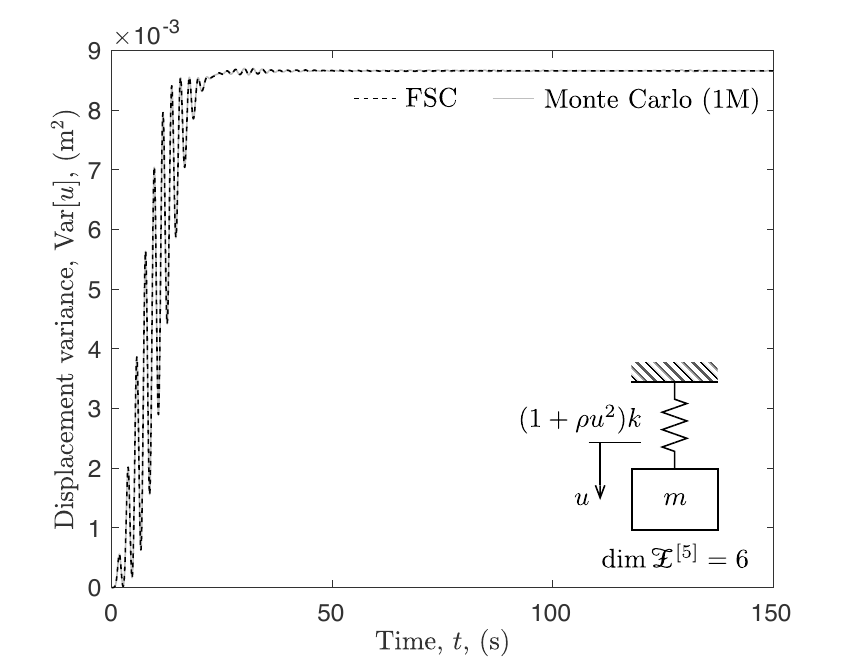}
\caption{Variance}
\label{fig1SDOF3_BetaUniform_FSC_Disp_Var_6}
\end{subfigure}
\caption{Evolution of $\mathbf{E}[u]$ and $\mathrm{Var}[u]$ for the case when the $p$-discretization level of RFS is $\mathscr{Z}^{[5]}$ and $\mu\sim\mathrm{Beta}\otimes\mathrm{Uniform}$}
\label{fig1SDOF3_BetaUniform_FSC_Disp_6}
\end{figure}

\begin{figure}
\centering
\begin{subfigure}[b]{0.495\textwidth}
\includegraphics[width=\textwidth]{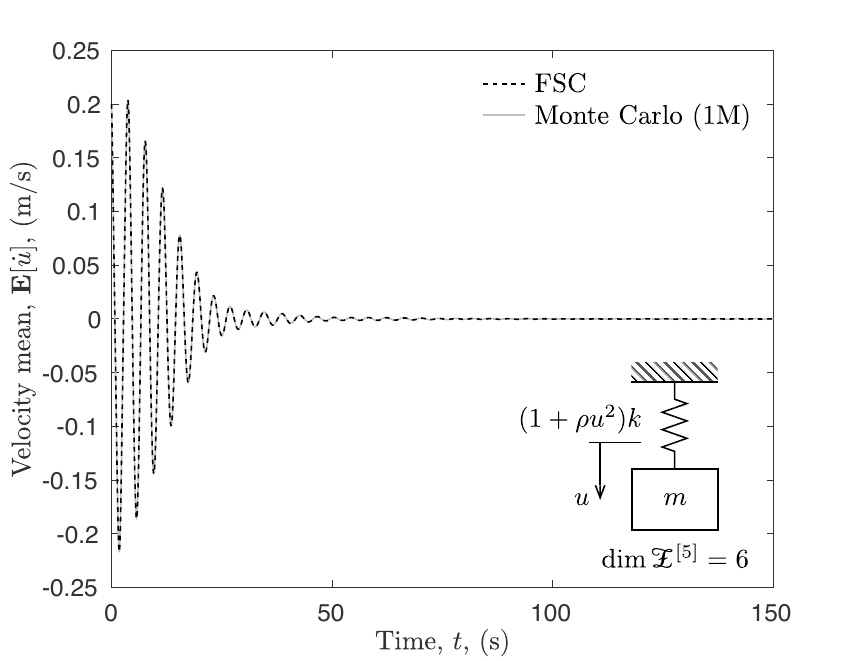}
\caption{Mean}
\label{fig1SDOF3_BetaUniform_FSC_Vel_Mean_6}
\end{subfigure}\hfill
\begin{subfigure}[b]{0.495\textwidth}
\includegraphics[width=\textwidth]{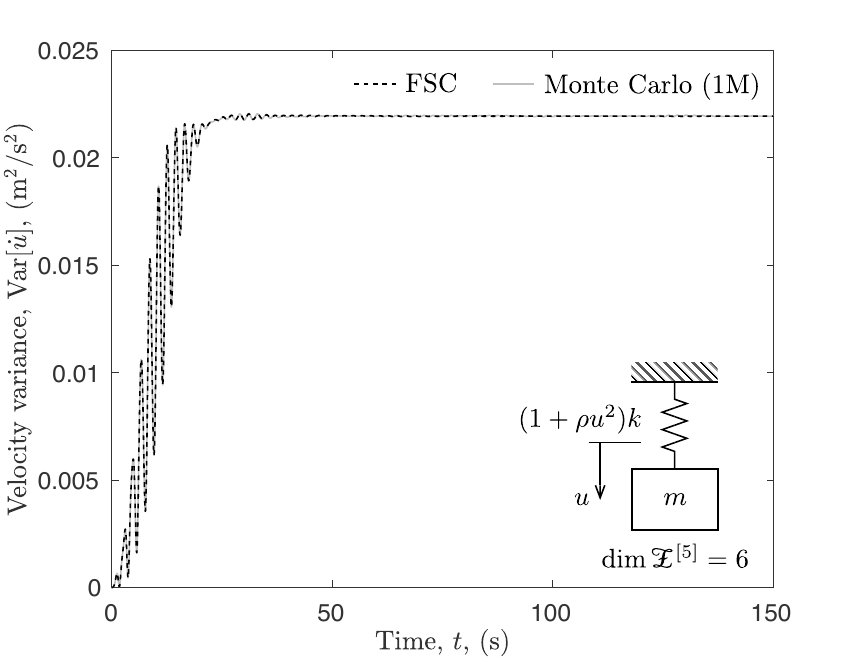}
\caption{Variance}
\label{fig1SDOF3_BetaUniform_FSC_Vel_Var_6}
\end{subfigure}
\caption{Evolution of $\mathbf{E}[\dot{u}]$ and $\mathrm{Var}[\dot{u}]$ for the case when the $p$-discretization level of RFS is $\mathscr{Z}^{[5]}$ and $\mu\sim\mathrm{Beta}\otimes\mathrm{Uniform}$}
\label{fig1SDOF3_BetaUniform_FSC_Vel_6}
\end{figure}

\begin{figure}
\centering
\begin{subfigure}[b]{0.495\textwidth}
\includegraphics[width=\textwidth]{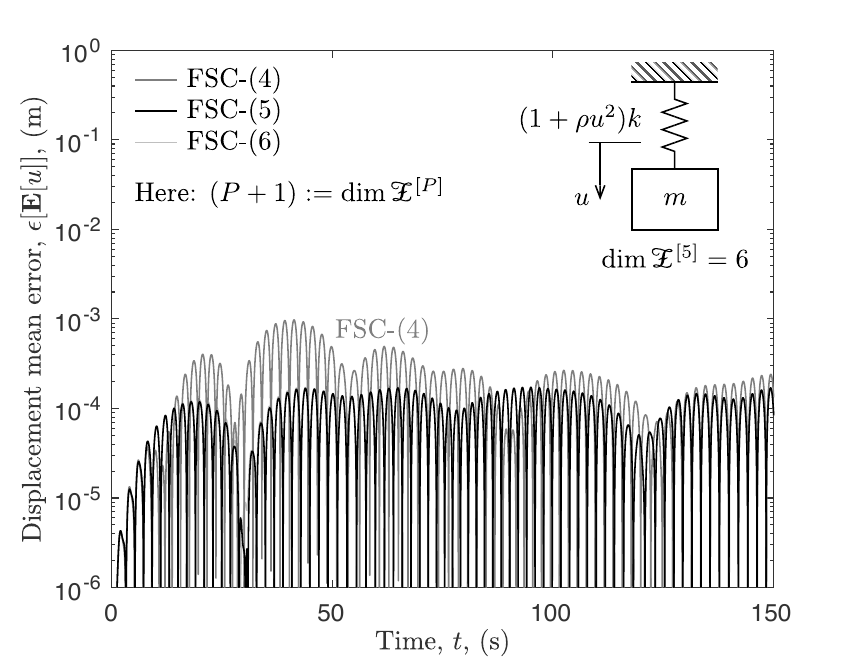}
\caption{Mean error}
\label{fig1SDOF3_BetaUniform_FSC_Disp_Mean_Error}
\end{subfigure}\hfill
\begin{subfigure}[b]{0.495\textwidth}
\includegraphics[width=\textwidth]{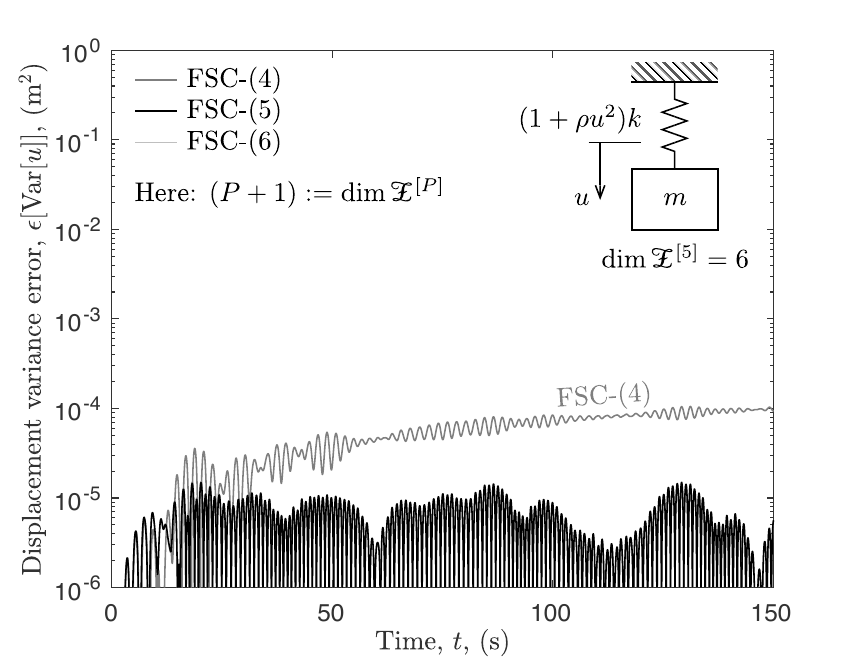}
\caption{Variance error}
\label{fig1SDOF3_BetaUniform_FSC_Disp_Var_Error}
\end{subfigure}
\caption{Local error evolution of $\mathbf{E}[u]$ and $\mathrm{Var}[u]$ for different $p$-discretization levels of RFS with respect to the 1-million Monte Carlo simulation ($\mu\sim\mathrm{Beta}\otimes\mathrm{Uniform}$)}
\label{fig1SDOF3_BetaUniform_FSC_Disp_Error}
\end{figure}

\begin{figure}
\centering
\begin{subfigure}[b]{0.495\textwidth}
\includegraphics[width=\textwidth]{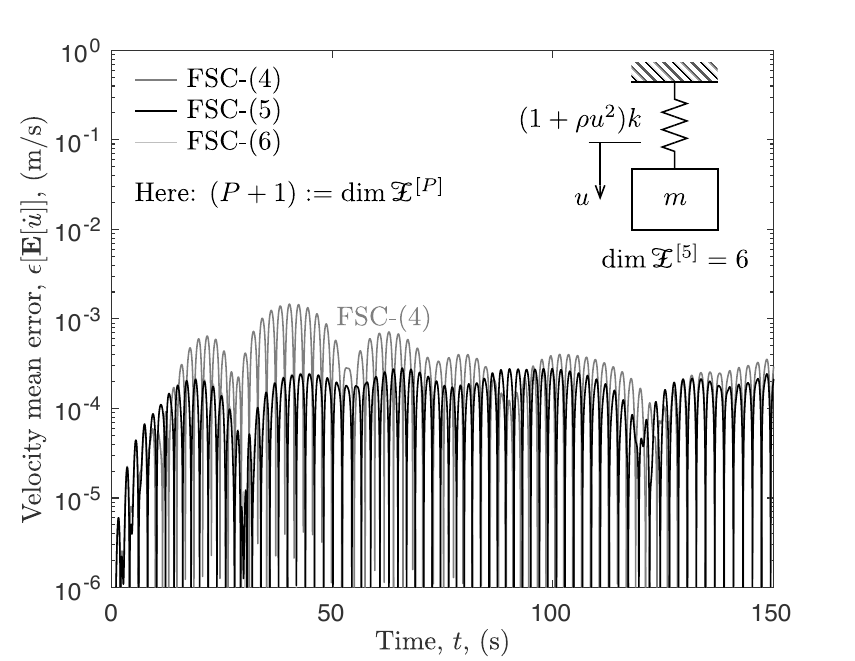}
\caption{Mean error}
\label{fig1SDOF3_BetaUniform_FSC_Vel_Mean_Error}
\end{subfigure}\hfill
\begin{subfigure}[b]{0.495\textwidth}
\includegraphics[width=\textwidth]{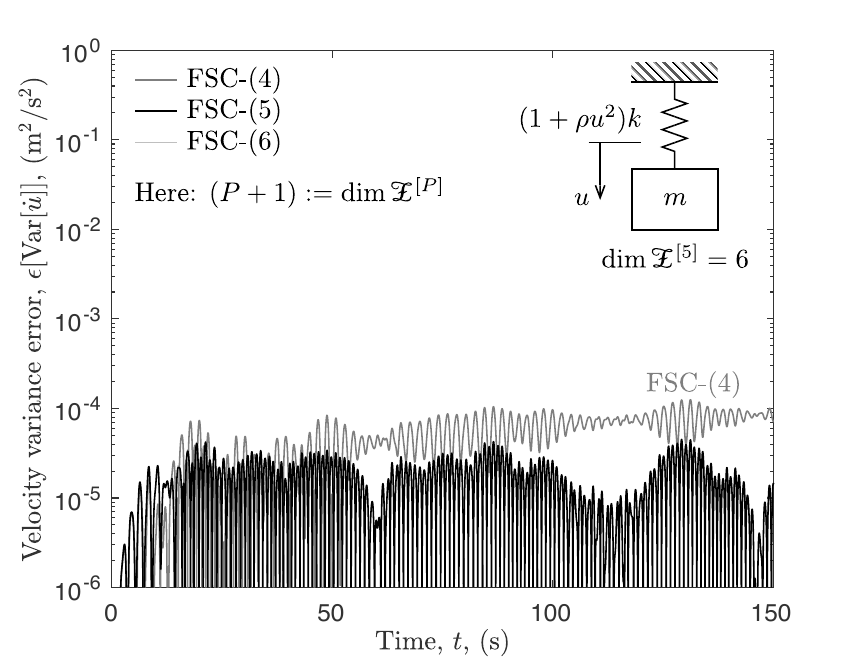}
\caption{Variance error}
\label{fig1SDOF3_BetaUniform_FSC_Vel_Var_Error}
\end{subfigure}
\caption{Local error evolution of $\mathbf{E}[\dot{u}]$ and $\mathrm{Var}[\dot{u}]$ for different $p$-discretization levels of RFS with respect to the 1-million Monte Carlo simulation ($\mu\sim\mathrm{Beta}\otimes\mathrm{Uniform}$)}
\label{fig1SDOF3_BetaUniform_FSC_Vel_Error}
\end{figure}

Since a closed-form solution for \eqref{eq1NumRes300} does not exist, we use the standard Monte Carlo method described in Appendix \ref{appsec1OveStaMonCarMet} to compare the accuracy of the FSC results against it.
To this end, one million realizations are randomly sampled from the random domain to conduct the Monte Carlo simulation.
The evolution of the mean and variance of the system's displacement is depicted in Figs.~\ref{fig1SDOF3_BetaUniform_FSC_Disp_6} and \ref{fig1SDOF3_BetaUniform_FSC_Vel_6}.
These figures show that by using only 5 basis vectors, the FSC results can effectively reproduce the Monte Carlo results.
In fact, Figs.~\ref{fig1SDOF3_BetaUniform_FSC_Disp_Error} and \ref{fig1SDOF3_BetaUniform_FSC_Vel_Error} further indicate that when $P=4,5$, the FSC results are 4-order-of-magnitude accurate for the mean and about 5-order-of-magnitude accurate for the variance.
However, this is not the case for $P=3$, which overall is one order of magnitude less accurate and displays a nearly linear drift after 50 s.
Moreover, it can be seen that the FSC results with $P=4$ and $P=5$ are indistinguishable from each other, chiefly because the one-million Monte Carlo simulation used as the reference solution is an approximate version of the solution.
This explains why the accuracy of the results did not increase when $P$ was set equal to 5.
Therefore, comparable results are achievable for this nonlinear problem if FSC is run with $P=4$.

\section{Application to structural dynamics}\label{sec1AppStrDyn}

\begin{figure}
\centering
\begin{subfigure}[b]{0.495\textwidth}
\includegraphics[width=\textwidth]{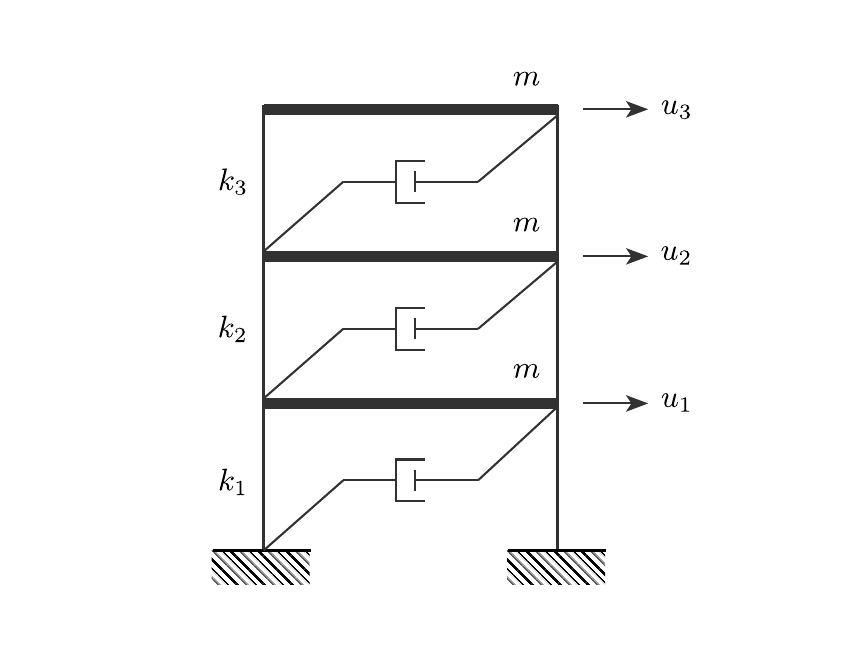}
\caption{Surrogate model of a 3-story building for lateral-load analysis in one direction}
\label{fig1threeStoryBuilding}
\end{subfigure}\hfill
\begin{subfigure}[b]{0.495\textwidth}
\includegraphics[width=\textwidth]{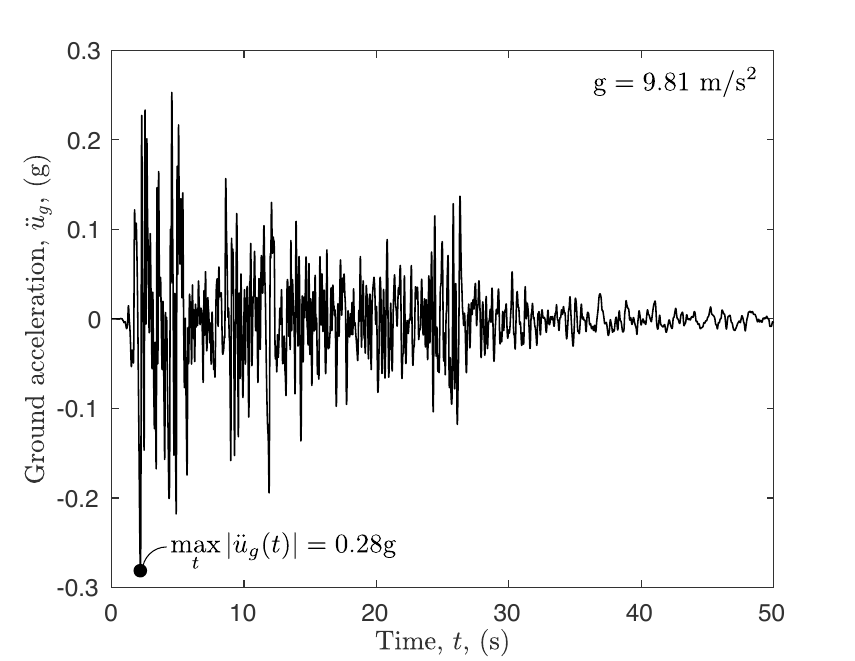}
\caption{Ground acceleration of \emph{1940 El Centro Earthquake}}
\label{fig1elCentroEarthquake1940}
\end{subfigure}
\caption{Surrogate model of a 3-story building to investigate its lateral behavior under an earthquake scenario}
\label{fig1MDOF3model}
\end{figure}

In structural dynamics, real-life systems are commonly modeled as multiple-degree-of-freedom systems.
In order to demonstrate how FSC can be utilized in a more general setting, in this section we quantify the response uncertainties of a 3-story building (Fig.~\ref{fig1threeStoryBuilding}) excited by the effects of a ground motion.
The ground motion is taken here to be one of the ground accelerations recorded from the \emph{1940 El Centro Earthquake}\footnote{This ground acceleration was obtained from the \emph{PEER Ground Motion Database} \cite{ancheta2014nga}.\\ Website: \url{https://ngawest2.berkeley.edu}.\\ Event's name: {\tt Imperial Valley-02}.
Station's name: {\tt El Centro Array \#9}.
File's name: {\tt RSN6\_IMPVALL.I\_I-ELC180.AT2}.}
event.
A plot of this ground motion is depicted in Fig.~\ref{fig1elCentroEarthquake1940} for sake of reference.

The governing differential equation of motion for this system is:
\begin{equation}\label{eq1AppStrDyn1000}
\mathbf{M}\mathbf{\ddot{u}}+\mathbf{C}\mathbf{\dot{u}}+\mathbf{K}\mathbf{u}=-\mathbf{M}\boldsymbol{\iota}\ddot{u}_g\quad(=:\mathbf{p}),
\end{equation}
where $\mathbf{M}\in\mathrm{L}(\mathbb{R}^3,\mathbb{R}^3)$ is the mass matrix, $\mathbf{C}:\Xi\to\mathrm{L}(\mathbb{R}^3,\mathbb{R}^3)$ is the damping matrix, $\mathbf{K}:\Xi\to\mathrm{L}(\mathbb{R}^3,\mathbb{R}^3)$ is the stiffness matrix, $\boldsymbol{\iota}\in\mathbb{R}^3$ is the influence vector, and $\ddot{u}_g:\mathfrak{T}\to\mathbb{R}$ is the ground acceleration characterized by a real-valued function of time.
The vectors $\mathbf{u},\mathbf{\dot{u}}:=\partial_t\mathbf{u},\mathbf{\ddot{u}}:=\partial_t^2\mathbf{u}:\mathfrak{T}\times\Xi\to\mathbb{R}^3$ represent, respectively, the displacement, the velocity and the acceleration of the system, where $\mathbf{u}^T=[u_1,u_2,u_3]$ is the unknown vector sought, and $u_3$ denotes the roof displacement of the 3-story building (our response of interest here).

The parameters of this system are defined as
\begin{equation*}
\mathbf{M}=m\begin{bmatrix}
1 & & \\
& 1 & \\
& & 1
\end{bmatrix},\quad
\mathbf{K}(\xi)=\begin{bmatrix}
k_1(\xi)+k_2(\xi) & -k_2(\xi) & \\
-k_2(\xi) & k_2(\xi)+k_3(\xi) & -k_3(\xi) \\
& -k_3(\xi) & k_3(\xi)
\end{bmatrix},
\end{equation*}
and $\mathbf{C}(\xi)=\alpha(\xi)\,\mathbf{M}+\beta(\xi)\,\mathbf{K}(\xi)$, where $m=500$ Mg, $k_1(\xi)=\xi^1\sim\mathrm{Beta}(2,5)$ in $[850\times10^3,1\,150\times10^3]$ kN/m, $k_2(\xi)=\xi^2\sim\mathrm{Beta}(2,5)$ in $[680\times10^3,920\times10^3]$ kN/m, $k_3(\xi)=\xi^3\sim\mathrm{Beta}(2,5)$ in $[680\times10^3,920\times10^3]$ kN/m, $\alpha(\xi)=\xi^4\sim\mathrm{Uniform}$ in $[0.4,0.7]$ s$^{-1}$, and $\beta(\xi)=\xi^5\sim\mathrm{Uniform}$ in $[0.4\times10^{-3},0.7\times10^{-3}]$ s.
Furthermore, $\boldsymbol{\iota}^T=[1,1,1]$, and $\ddot{u}_g$ is defined according to Fig.~\ref{fig1elCentroEarthquake1940}.
The random domain for this system is thus 5-dimensional:
\begin{equation*}
\Xi=\prod_{i=1}^5 \bar{\Xi}_i\equiv[850\times10^3,1\,150\times10^3]\times[680\times10^3,920\times10^3]^2\times[0.4,0.7]\times[0.4\times10^{-3},0.7\times10^{-3}],
\end{equation*}
and we assume that the initial state of the system is at rest, i.e.~$\mathbf{u}(0)=\mathbf{\dot{u}}(0)=\mathbf{0}$.
Note that when the expected values of $k_1$, $k_2$ and $k_3$ are utilized to define $\mathbf{K}$, the fundamental period of the system is approximately 0.33 s, which is consistent with a typical 3-story building found in practice featuring a damping ratio of about 2\% for the first two modal frequencies.

Following the formulation presented in Section \ref{sec1SolSpeApp}, the system of equations \eqref{eq1SolSpeApp1060} takes the form:
\begin{subequations}\label{eq1AppStrDyn1060}
\begin{align}
\mathbf{M}\indices{^i_j}\mathbf{\ddot{u}}^j+\mathbf{C}\indices{^i_j}\mathbf{\dot{u}}^j+\mathbf{K}\indices{^i_j}\mathbf{u}^j=\mathbf{p}^i &\qquad\text{on $\mathfrak{T}=[0,50]$ s}\label{eq1AppStrDyn1060a}\\
\big\{\mathbf{u}^i(0)=\mathbf{0},\,\mathbf{\dot{u}}^i(0)=\mathbf{0}\big\}&\qquad\text{on $\{0\}$},\label{eq1AppStrDyn1060b}
\end{align}
\end{subequations}
where a summation sign is implied over the repeated index $j$, $\mathbf{M}\indices{^i_j}=\mathbf{M}\delta\indices{^i_j}$, $\mathbf{p}^i(t)=-\mathbf{M}\boldsymbol{\iota}\,\ddot{u}_g(t)\,\delta\indices{^i_0}$,
\begin{multline*}
\mathbf{C}\indices{^i_j}=
m\begin{bmatrix}
1 & & \\
 & 1 & \\
 & & 1
\end{bmatrix}\frac{\langle\Psi_i,\alpha\Psi_j\rangle}{\langle\Psi_i,\Psi_i\rangle}
+
\begin{bmatrix}
1 & & \\
 & 0 & \\
 & & 0
\end{bmatrix}\frac{\langle\Psi_i,\beta k_1\Psi_j\rangle}{\langle\Psi_i,\Psi_i\rangle}\\[1.5ex]
+
\begin{bmatrix}
1 & -1 & \\
-1 & 1 &  \\
 &  & 0
\end{bmatrix}\frac{\langle\Psi_i,\beta k_2\Psi_j\rangle}{\langle\Psi_i,\Psi_i\rangle}
+
\begin{bmatrix}
0 &  & \\
 & 1 & -1 \\
 & -1 & 1
\end{bmatrix}\frac{\langle\Psi_i,\beta k_3\Psi_j\rangle}{\langle\Psi_i,\Psi_i\rangle},
\end{multline*}
and
\begin{equation*}
\mathbf{K}\indices{^i_j}=
\begin{bmatrix}
1 & & \\
 & 0 & \\
 & & 0
\end{bmatrix}\frac{\langle\Psi_i,k_1\Psi_j\rangle}{\langle\Psi_i,\Psi_i\rangle}
+
\begin{bmatrix}
1 & -1 & \\
-1 & 1 &  \\
 &  & 0
\end{bmatrix}\frac{\langle\Psi_i,k_2\Psi_j\rangle}{\langle\Psi_i,\Psi_i\rangle}
+
\begin{bmatrix}
0 &  & \\
 & 1 & -1 \\
 & -1 & 1
\end{bmatrix}\frac{\langle\Psi_i,k_3\Psi_j\rangle}{\langle\Psi_i,\Psi_i\rangle}
\end{equation*}
with $i,j\in\{0,1,\ldots,P\}$, and $\delta\indices{^i_j}$ denoting the Kronecker delta.

For this problem, we take $P=9$ for the set of linearly independent functions $\{\Phi_{j.i}:=\hat{\varphi}^j(0,\mathbf{\hat{s}}(t_i,\cdot\,))\}_{j=1}^P$, where $\hat{\varphi}:\mathbb{R}\times\mathscr{Z}^9\to\mathscr{Z}^9$ is given by
\begin{equation*}
\hat{\varphi}(h,\mathbf{\hat{s}}(t_i,\cdot\,))=:\mathbf{\hat{s}}(t_i+h,\cdot\,)=(\mathbf{u}(t_i+h,\cdot\,),\mathbf{\dot{u}}(t_i+h,\cdot\,),\mathbf{\ddot{u}}(t_i+h,\cdot\,)),
\end{equation*}
and $\mathbf{\hat{s}}=(\mathbf{u},\mathbf{\dot{u}},\mathbf{\ddot{u}})\equiv(u_1,u_2,u_3,\dot{u}_1,\dot{u}_2,\dot{u}_3,\ddot{u}_1,\ddot{u}_2,\ddot{u}_3)$.

\begin{remark}
Notice that this $P$ is the smallest value we can choose for a problem featuring three degrees of freedom and whose governing stochastic differential equation is of second-order in time.
This is because for each degree of freedom, the smallest RFS that one can construct using the FSC method is one whose $P$ is equal to 3.
However, we emphasize that this issue is not particular to the FSC method. 
It is well recognized that the more degrees of freedom a dynamical system has, the more basis vectors needed to construct a suitable random function space for the system's state at any given time.
For example, in TD-gPC-based methods, this would be equivalent to perform a full tensor product between all the RFS's generated at each degree of freedom.	
\end{remark}

To integrate \eqref{eq1AppStrDyn1060} numerically, we employ the RK4 method with a time-step size of $\Delta t=0.01$ s (which is concordant with the sample frequency of the ground acceleration record).
To evaluate the inner products approximately, we use 15 Gaussian quadrature points on each random axis, resulting in $15^5=759\,375$ quadrature points distributed across the entire random domain.

\begin{figure}
\centering
\begin{subfigure}[b]{0.495\textwidth}
\includegraphics[width=\textwidth]{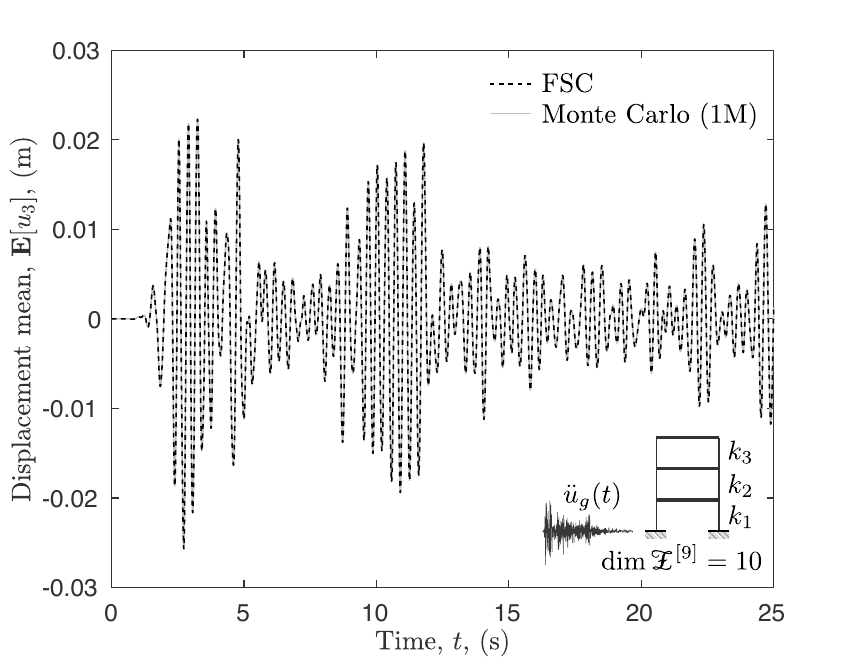}
\caption{Mean}
\label{fig1MDOF3_Uniform2_FSC_Disp_Mean_10}
\end{subfigure}\hfill
\begin{subfigure}[b]{0.495\textwidth}
\includegraphics[width=\textwidth]{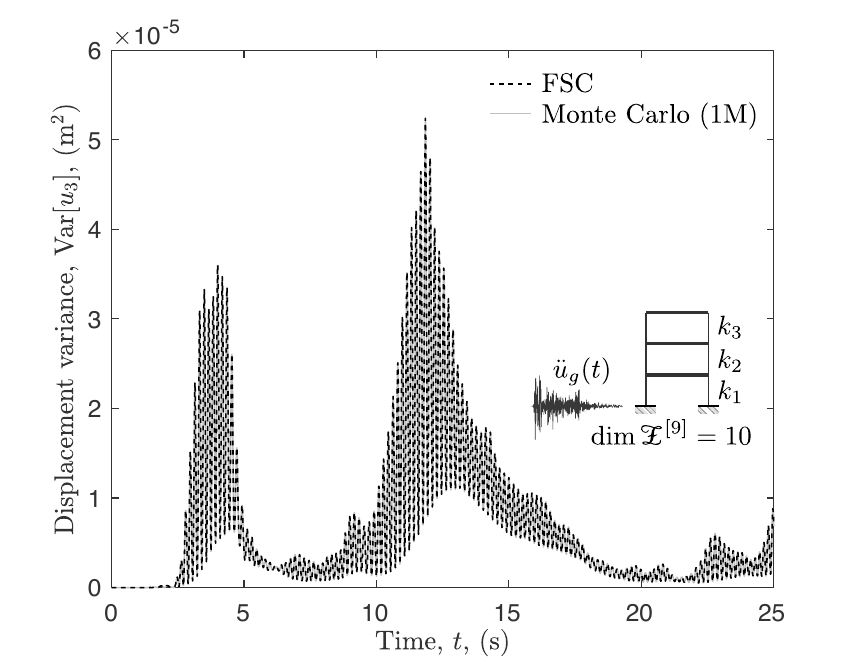}
\caption{Variance}
\label{fig1MDOF3_Uniform2_FSC_Disp_Var_10}
\end{subfigure}
\caption{Evolution of $\mathbf{E}[u_3]$ and $\mathrm{Var}[u_3]$ for the case when the $p$-discretization level of RFS is $\mathscr{Z}^{[9]}$ and $\mu\sim\mathrm{Beta}^{\otimes3}\otimes\mathrm{Uniform}^{\otimes2}$}
\label{fig1MDOF3_Uniform2_FSC_Disp_10}
\end{figure}

\begin{figure}
\centering
\begin{subfigure}[b]{0.495\textwidth}
\includegraphics[width=\textwidth]{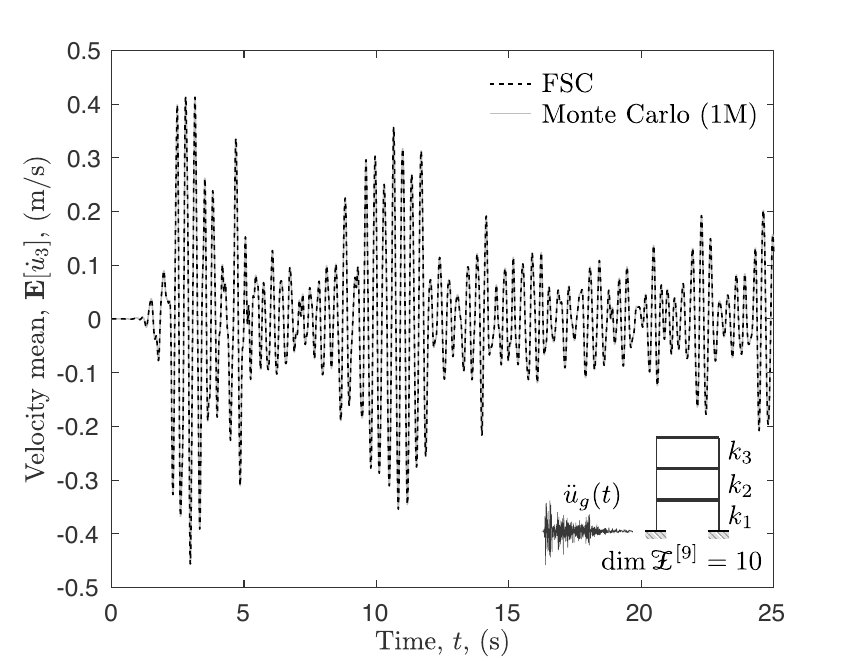}
\caption{Mean}
\label{fig1MDOF3_Beta3Uniform2_FSC_Vel_Mean_10}
\end{subfigure}\hfill
\begin{subfigure}[b]{0.495\textwidth}
\includegraphics[width=\textwidth]{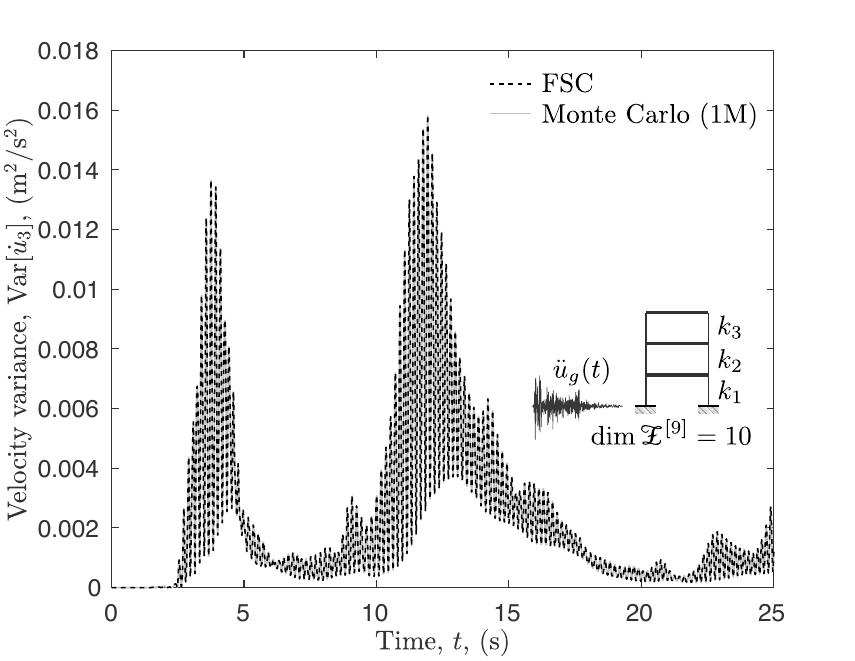}
\caption{Variance}
\label{fig1MDOF3_Beta3Uniform2_FSC_Vel_Var_10}
\end{subfigure}
\caption{Evolution of $\mathbf{E}[\dot{u}_3]$ and $\mathrm{Var}[\dot{u}_3]$ for the case when the $p$-discretization level of RFS is $\mathscr{Z}^{[9]}$ and $\mu\sim\mathrm{Beta}^{\otimes3}\otimes\mathrm{Uniform}^{\otimes2}$}
\label{fig1MDOF3_Beta3Uniform2_FSC_Vel_10}
\end{figure}

The results in Figs.~\ref{fig1MDOF3_Uniform2_FSC_Disp_10} and \ref{fig1MDOF3_Beta3Uniform2_FSC_Vel_10} depict the solutions obtained from employing FSC and the standard Monte Carlo method (as described in Appendix \ref{appsec1OveStaMonCarMet}) to quantify the uncertainties of the response.
One million realizations were randomly sampled from the random domain to conduct the Monte Carlo simulation.
Once again, it is apparent that FSC is able to capture the system's uncertainties with high fidelity and its solution is indistinguishable from that of the Monte Carlo method.
For clarity, only the first 25 s of the solution are presented in Figs.~\ref{fig1MDOF3_Uniform2_FSC_Disp_10} and \ref{fig1MDOF3_Beta3Uniform2_FSC_Vel_10}, however, the conclusion drawn above applies to the last 25 s of the solution as well.

\section{Conclusion}

A novel numerical method, called the \emph{flow-driven spectral chaos} (FSC) method, is presented for capturing uncertainties in structural dynamics using the spectral approach.
The FSC method uses the concept of enriched stochastic flow maps to track the evolution of the system's state efficiently in an augmented random phase space.
The method is not only computationally more efficient than the TD-gPC approach but also easy to implement, since the flow map that we use in the scheme is nothing but the time derivatives of the solution up to a specific order.
Moreover, since the random basis is defined with these time derivatives, the number of basis vectors required to characterize the stochastic part of the solution space does not depend upon the dimensionality of the probability space.
This remarkable property opens up the possibility of investigating systems with high-dimensional probability spaces at low computational cost---an issue that has plagued the spectral approach since the introduction of the PC method.

The three numerical examples presented in Section \ref{sec1NumRes} show that the FSC scheme is able to capture the response of the system with high accuracy using a small number of basis vectors and at a relatively low computational cost.
The illustrative problem described in Section \ref{sec1AppStrDyn} also shows that the FSC method can be readily applied to real-world structures involving multiple degrees of freedom.
As a result, the FSC method has the potential to be used in the context of large-scale structural engineering problems to quantify the uncertainties of long-time response with high fidelity.

\begin{appendices}
\section{Random bases for illustrative examples}\label{appsec1RanBasTimInt}

Table \ref{apptab1RanBasTimInt1000} presents the non-orthogonalized version of the random bases that we use in this manuscript to solve the examples described in Section \ref{sec1NumRes}.

\begin{table}
\centering\small
\caption{Non-orthogonalized version of the random bases used in Section \ref{sec1NumRes}}
\label{apptab1RanBasTimInt1000}
\begin{tabular}{@{}p{\textwidth}@{}}
\toprule
Single-degree-of-freedom system under {\bf free vibration} (Section \ref{sec1NumRes10})$^*$:
\begin{gather*}
\Phi_{0.i}(\xi):=1\\
\Phi_{1.i}(\xi):=\hat{\varphi}^1(M)(0,\hat{s}(t_i,\xi))= u_{.i}(t_i,\xi)=u_{.i-1}(t_i,\xi)\\
\Phi_{2.i}(\xi):=\hat{\varphi}^2(M)(0,\hat{s}(t_i,\xi))= \dot{u}_{.i}(t_i,\xi)=\dot{u}_{.i-1}(t_i,\xi)\\
\Phi_{3.i}(\xi):=\hat{\varphi}^3(M)(0,\hat{s}(t_i,\xi))=\partial_t^2u_{.i}(t_i,\xi)=-\frac{k(\xi)}{m}\,u_{.i}(t_i,\xi)\\
\Phi_{4.i}(\xi):=\hat{\varphi}^4(M)(0,\hat{s}(t_i,\xi))=\partial_t^3u_{.i}(t_i,\xi)=-\frac{k(\xi)}{m}\,\dot{u}_{.i}(t_i,\xi)\\
\Phi_{5.i}(\xi):=\hat{\varphi}^5(M)(0,\hat{s}(t_i,\xi))=\partial_t^4u_{.i}(t_i,\xi)=-\frac{k(\xi)}{m}\,\partial_t^2u_{.i}(t_i,\xi)\\
\Phi_{6.i}(\xi):=\hat{\varphi}^6(M)(0,\hat{s}(t_i,\xi))=\partial_t^5u_{.i}(t_i,\xi)=-\frac{k(\xi)}{m}\,\partial_t^3u_{.i}(t_i,\xi)\\
\vdots\quad\text{(until $P=M+2$ if needed)}
\end{gather*}\\[-3ex]
\midrule
Single-degree-of-freedom system under {\bf forced vibration} (Section \ref{sec1NumRes20})$^*$:
\begin{gather*}
\Phi_{0.i}(\xi):=1\\
\Phi_{1.i}(\xi):=\hat{\varphi}^1(M)(0,\hat{s}(t_i,\xi))=u_{.i}(t_i,\xi)=u_{.i-1}(t_i,\xi)\\
\Phi_{2.i}(\xi):=\hat{\varphi}^2(M)(0,\hat{s}(t_i,\xi))=\dot{u}_{.i}(t_i,\xi)=\dot{u}_{.i-1}(t_i,\xi)\\
\Phi_{3.i}(\xi):=\hat{\varphi}^3(M)(0,\hat{s}(t_i,\xi))=\partial_t^2u_{.i}(t_i,\xi)=\frac{1}{m}\big(q(\xi)\,\sin(t_i)-k(\xi)\,u_{.i}(t_i,\xi)\big)\\
\Phi_{4.i}(\xi):=\hat{\varphi}^4(M)(0,\hat{s}(t_i,\xi))=\partial_t^3u_{.i}(t_i,\xi)=\frac{1}{m}\big(q(\xi)\,\cos(t_i)-k(\xi)\,\dot{u}_{.i}(t_i,\xi)\big)\\
\Phi_{5.i}(\xi):=\hat{\varphi}^5(M)(0,\hat{s}(t_i,\xi))=\partial_t^4u_{.i}(t_i,\xi)=\frac{1}{m}\big({-q(\xi)}\,\sin(t_i)-k(\xi)\,\partial_t^2u_{.i}(t_i,\xi)\big)\\
\Phi_{6.i}(\xi):=\hat{\varphi}^6(M)(0,\hat{s}(t_i,\xi))=\partial_t^5u_{.i}(t_i,\xi)=\frac{1}{m}\big({-q(\xi)}\,\cos(t_i)-k(\xi)\,\partial_t^3u_{.i}(t_i,\xi)\big)\\
\vdots\quad\text{(until $P=M+2$ if needed)}
\end{gather*}\\[-3ex]
\midrule
{\bf Nonlinear} single-degree-of-freedom system under {\bf free vibration} (Section \ref{sec1NumRes30})$^*$:
\begin{gather*}
\Phi_{0.i}(\xi):=1\\
\Phi_{1.i}(\xi):=\hat{\varphi}^1(M)(0,\hat{s}(t_i,\xi))=u_{.i}(t_i,\xi)=u_{.i-1}(t_i,\xi)\\
\Phi_{2.i}(\xi):=\hat{\varphi}^2(M)(0,\hat{s}(t_i,\xi))=\dot{u}_{.i}(t_i,\xi)=\dot{u}_{.i-1}(t_i,\xi)\\
\Phi_{3.i}(\xi):=\hat{\varphi}^3(M)(0,\hat{s}(t_i,\xi))=\partial_t^2u_{.i}(t_i,\xi)=-\frac{k(\xi)}{m}(1+\rho(\xi)\,u^2{_{.i}}(t_i,\xi))\,u_{.i}(t_i,\xi)\\
\Phi_{4.i}(\xi):=\hat{\varphi}^4(M)(0,\hat{s}(t_i,\xi))=\partial_t^3u_{.i}(t_i,\xi)=-\frac{k(\xi)}{m}(1+3\,\rho(\xi)\,u^2{_{.i}}(t_i,\xi))\,\dot{u}_{.i}(t_i,\xi)\\
\Phi_{5.i}(\xi):=\hat{\varphi}^5(M)(0,\hat{s}(t_i,\xi))=\partial_t^4u_{.i}(t_i,\xi)=-\frac{k(\xi)}{m}(1+3\,\rho(\xi)\,u^2{_{.i}}(t_i,\xi))\,\partial_t^2 u_{.i}(t_i,\xi)-\frac{6\,\rho(\xi)\,k(\xi)}{m}\,\dot{u}^2{_{.i}}(t_i,\xi)\,u_{.i}(t_i,\xi)\\
\vdots\quad\text{(until $P=M+2$ if needed)}
\end{gather*}\\[-3ex]
\midrule
{\footnotesize $^*$The random basis is defined over the region $\mathfrak{R}_i=\mathrm{cl}(\mathfrak{T}_i)\times\Xi$.}\\
\bottomrule
\end{tabular}
\end{table}

\section{Undamped single-degree-of-freedom system under free vibration}\label{appsec1SDOF}

The objective of this section is to provide the exact response expressions of an undamped single-degree-of-freedom system subjected to free vibration for the case when the stiffness is assumed uniformly distributed with parameters $k_a$ and $k_b$.
The problem is stated formally as follows.

\paragraph{Problem statement} Consider a stochastic, undamped single-degree-of-freedom system with mass $m\in\mathbb{R}^+$, and stiffness $k:\Xi\to\mathbb{R}^+$ given by $k(\xi)=\xi\sim\mathrm{Uniform}$ in $[k_a,k_b]$, subjected to free vibration. Note that: $k_b>k_a>0$.

The \emph{first problem} is to find the displacement of the system $u:\mathfrak{T}\times\Xi\to\mathbb{R}$ in $\mathscr{U}$, such that:
\begin{subequations}\label{appeq1SDOF1000}
\begin{align}
m\ddot{u}+ku=0&\qquad\text{on $\mathfrak{T}\times\Xi$}\label{appeq1SDOF1000a}\\
\big\{u(0,\cdot\,)=\mathscr{u},\,\dot{u}(0,\cdot\,)=\mathscr{v}\big\}&\qquad\text{on $\{0\}\times\Xi$}\label{appeq1SDOF1000b}
\end{align}
\end{subequations}
where $\mathscr{u},\mathscr{v}\in\mathbb{R}$, and $\dot{u}:=\partial_t u$ and $\ddot{u}:=\partial_t^2 u$ are the velocity and acceleration of the system, respectively.

The \emph{second problem} is to find the expectation and variance of $u$, $\dot{u}$ and $\ddot{u}$ as a function of time.

\paragraph{Exact solution} The solution of \eqref{appeq1SDOF1000} is well-known \cite{chopra2012dynamics,humar2012dynamics}, and it is given by
\begin{equation}\label{appSDOF3000}
u(t,\xi)=\mathscr{u}\cos\!\big((\omega\circ k)(\xi)\,t\big)+\frac{\mathscr{v}}{(\omega\circ k)(\xi)}\sin\!\big((\omega\circ k)(\xi)\,t\big),
\end{equation}
where the natural circular frequency of the system, $\omega:\mathbb{R}^+\to\mathbb{R}^+$, is defined by
\begin{equation*}
\omega(k)=\sqrt{\frac{k}{m}}.
\end{equation*}

Then, the following exact expressions can be derived for the expectation and variance of $u$, $\dot{u}$ and $\ddot{u}$.

\paragraph{Exact expectation} The expectation of $u$, $\dot{u}$ and $\ddot{u}$ are given by:
\begin{subequations}\label{appeq1SDOF5000}
\begin{align}
\mathbf{E}[u](t)&=\kappa(t)\,(\tau_u(t,k_b)-\tau_u(t,k_a))\\
\mathbf{E}[\dot{u}](t)&=\kappa(t)\,(\tau_v(t,k_b)-\tau_v(t,k_a))\\
\mathbf{E}[\ddot{u}](t)&=\kappa(t)\,(\tau_a(t,k_b)-\tau_a(t,k_a)),
\end{align}
\end{subequations}
where $\kappa:\mathfrak{T}\to\mathbb{R}^+$ is given by $\kappa(t)=2m/((k_b-k_a)t^2)$, and $\tau_u,\tau_v,\tau_a:\mathfrak{T}\times\mathbb{R}^+\to\mathbb{R}$ are defined in Table \ref{apptab1SDOF1000}.

\begin{table}
\centering\small
\caption{Definition of $\tau$-functions for a single-degree-of-freedom system subjected to free vibration with $k\sim\mathrm{Uniform}$ in $[k_a,k_b]$}
\label{apptab1SDOF1000}
\begin{tabular}{@{}p{\textwidth}@{}}
\toprule\vspace{-2ex}
\begin{equation*}
\tau_u(t,k)=\big\{\omega(k)\,t\sin(\omega(k)\,t)+\cos(\omega(k)\,t)\big\}\,\mathscr{u}
-\cos(\omega(k)\,t)\,\mathscr{v}t
\end{equation*}\\[-3ex]
\midrule\vspace{-2ex}
\begin{equation*}
\tau_v(t,k)=-\big\{2\,\omega(k)\,t\sin(\omega(k)\,t)+(2-\omega^2(k)\,t^2)\cos(\omega(k)\,t)\big\}\,\mathscr{u}t^{-1}
+\big\{\omega(k)\,t\sin(\omega(k)\,t)+\cos(\omega(k)\,t)\big\}\,\mathscr{v}	
\end{equation*}\\[-3ex]
\midrule\vspace{-5ex}
\begin{multline*}
\tau_a(t,k)=-\big\{(\omega^3(k)\,t^3-6\,\omega(k)\,t)\sin(\omega(k)\,t)+3(\omega^2(k)\,t^2-2)\cos(\omega(k)\,t)\big\}\,\mathscr{u}t^{-2}\\
-\big\{2\,\omega(k)\,t\sin(\omega(k)\,t)+(2-\omega^2(k)\,t^2)\cos(\omega(k)\,t)\big\}\,\mathscr{v}t^{-1}
\end{multline*}\\[-4ex]
\bottomrule
\end{tabular}
\end{table}

A closer look at the above expressions indicates that in the long term the absolute mean of the response is dominated by $\kappa=\kappa(t)$ which is a function that tends to zero as time goes to infinity.
For this reason,
\begin{equation}\label{appeq1SDOF5500}
\mathbf{E}[u],\mathbf{E}[\dot{u}],\mathbf{E}[\ddot{u}]\to0\quad\text{as}\quad t\to\infty.
\end{equation}

\paragraph{Exact variance} The variance of $u$, $\dot{u}$ and $\ddot{u}$ are given by:
\begin{subequations}\label{appeq1SDOF6000}
\begin{align}
\mathrm{Var}[u](t)&=\kappa(t)\,\big(\varrho_u(t,k_b)-\varrho_u(t,k_a)\big)-\mathbf{E}[u]^2(t)\\
\mathrm{Var}[\dot{u}](t)&=\kappa(t)\,\big(\varrho_v(t,k_b)-\varrho_v(t,k_a)\big)-\mathbf{E}[\dot{u}]^2(t)\\
\mathrm{Var}[\ddot{u}](t)&=\kappa(t)\,\big(\varrho_a(t,k_b)-\varrho_a(t,k_a)\big)-\mathbf{E}[\ddot{u}]^2(t),
\end{align}
\end{subequations}
where $\varrho_u,\varrho_v,\varrho_a:\mathfrak{T}\times\mathbb{R}^+\to\mathbb{R}$ are defined in Table \ref{apptab1SDOF2000}. In the expression for $\varrho_u$, the function $\mathrm{Ci}:\mathbb{R}^+\to\mathbb{R}$ is the cosine integral given by
\begin{equation*}
\mathrm{Ci}(x)=-\int_x^\infty\frac{\cos y}{y}\,\mathrm{d}y.
\end{equation*}

\begin{table}
\centering\small
\caption{Definition of $\varrho$-functions for a single-degree-of-freedom system subjected to free vibration with $k\sim\mathrm{Uniform}$ in $[k_a,k_b]$}
\label{apptab1SDOF2000}
\begin{tabular}{@{}p{\textwidth}@{}}
\toprule\vspace{-5ex}
\begin{multline*}
\varrho_u(t,k)=-\tfrac{1}{4}\big\{\!\sin^2(\omega(k)\,t)-\omega(k)\,t\sin(2\,\omega(k)\,t)-\omega^2(k)\,t^2\big\}\,\mathscr{u}^2
-\tfrac{1}{2}\cos(2\,\omega(k)\,t)\,\mathscr{u}\mathscr{v}t\\
+\tfrac{1}{4}\big\{\!\ln(k)-2\,\mathrm{Ci}(2\,\omega(k)\,t)\big\}\,\mathscr{v}^2t^2
\end{multline*}\\[-4ex]
\midrule\vspace{-5ex}
\begin{multline*}
\varrho_v(t,k)=\tfrac{1}{16}\big\{2\,(3\,\omega(k)\,t-2\,\omega^3(k)\,t^3)\sin(2\,\omega(k)\,t)+3\,(1-2\,\omega^2(k)\,t^2)\cos(2\,\omega(k)\,t)+2\,\omega^4(k)\,t^4\big\}\,\mathscr{u}^2 t^{-2}\\
-\tfrac{1}{4}\big\{2\,\omega(k)\,t\sin(2\,\omega(k)\,t)+(1-2\,\omega^2(k)\,t^2)\cos(2\,\omega(k)\,t)\big\}\,\mathscr{u}\mathscr{v}t^{-1}\\
-\tfrac{1}{4}\big\{\!\sin^2(\omega(k)\,t)-\omega(k)\,t\sin(2\,\omega(k)\,t)-\omega^2(k)\,t^2\big\}\,\mathscr{v}^2
\end{multline*}\\[-4ex]
\midrule\vspace{-5ex}
\begin{multline*}
\varrho_a(t,k)=\tfrac{1}{48}\big\{6\,(2\,\omega^5(k)\,t^5-10\,\omega^3(k)\,t^3+15\,\omega(k)\,t)\sin(2\,\omega(k)\,t)\\
+15\,(2\,\omega^4(k)\,t^4-6\,\omega^2(k)\,t^2+3)\cos(2\,\omega(k)\,t)+4\,\omega^6(k)\,t^6\big\}\,\mathscr{u}^2t^{-4}\\
+\tfrac{1}{8}\big\{4\,(2\,\omega^3(k)\,t^3-3\,\omega(k)\,t)\sin(2\,\omega(k)\,t)-2\,(2\,\omega^4(k)\,t^4-6\,\omega^2(k)\,t^2+3)\cos(2\,\omega(k)\,t)\big\}\,\mathscr{u}\mathscr{v}t^{-3}\\
-\tfrac{1}{16}\big\{2\,(2\,\omega^3(k)\,t^3-3\,\omega(k)\,t)\sin(2\,\omega(k)\,t)+3\,(2\,\omega^2(k)\,t^2-1)\cos(2\,\omega(k)\,t)-2\,\omega^4(k)\,t^4\big\}\,\mathscr{v}^2t^{-2}
\end{multline*}\\[-4ex]
\bottomrule
\end{tabular}
\end{table}

We see that not only the variance of the response is bounded for all $t\in\mathfrak{T}$, but also as $t$ goes to infinity:
\begin{subequations}\label{appeq1SDOF7000}
\begin{gather}
\lim_{t\to\infty}\mathrm{Var}[u](t)=\frac{1}{2}\mathscr{u}^2+\frac{1}{2}\ln\!\left(\frac{k_b}{k_a}\right)\!\left(\frac{m}{k_b-k_a}\right)\!\mathscr{v}^2\\
\lim_{t\to\infty}\mathrm{Var}[\dot{u}](t)=\frac{1}{4}\!\left(\frac{k_a+k_b}{m}\right)\!\mathscr{u}^2+\frac{1}{2}\mathscr{v}^2\\
\lim_{t\to\infty}\mathrm{Var}[\ddot{u}](t)=\frac{1}{6}\!\left(\frac{k_a^2+k_ak_b+k_b^2}{m^2}\right)\!\mathscr{u}^2+\frac{1}{4}\!\left(\frac{k_a+k_b}{m}\right)\!\mathscr{v}^2.
\end{gather}
\end{subequations}
Therefore, the underlying process $u$ is naturally of second-order.

\section{Overview of standard Monte Carlo method}\label{appsec1OveStaMonCarMet}

The Monte Carlo method is the most popular numerical technique used in stochastic modeling to quantify the effects of input uncertainty on system's outputs.
It is basically a `brute-force' method of attack that typically involves sampling a large number of realizations from the random space to estimate the statistics of the output.
It is well-known that when $N$ realizations are considered, the mean converges asymptotically as the square root of $N^{-1}$, and thus, it is remarkably independent of the dimensionality of the random space \cite{fishman2013monte}.
In this paper we use \emph{standard Monte Carlo} to validate the FSC method in Sections \ref{sec1NumRes30} and \ref{sec1AppStrDyn}.

The general procedure for conducting a standard Monte Carlo simulation is simple.
Consider the stochastic system given by \eqref{eq1ProSta1000star}:
\begin{equation}\label{eq1ProSta1000star2}
y(t,\xi)=\boldsymbol{\mathcal{M}}[u][x](t,\xi)\quad\text{subject to initial condition}\quad \boldsymbol{\mathcal{I}}[u](\xi).\tag{\ref{eq1ProSta1000star}*}
\end{equation}
Then:
\begin{enumerate}
\item Generate $N$ realizations of the $d$-tuple random variable $\xi$ in order to obtain the random set $\{\xi_i\}_{i=1}^N$.
These $N$ realizations are based on randomly sampling $N$ points from the random domain $\Xi$ according to the cumulative distribution function $F:\Xi\to[0,1]$ given by
\begin{equation}
F(\xi)=\prod_{j=1}^d\mu^j\big((-\infty,\xi^j]\big),\quad\text{or equivalently,}\quad
F(\xi)=\int_{-\infty}^{\xi^1}\cdots\int_{-\infty}^{\xi^d}\mathrm{d}\mu^1\cdots\mathrm{d}\mu^d.
\end{equation}
This way we can also obtain the input set $\{x(t,\xi_i)\}_{i=1}^N$ for reference purposes.
\item Solve \eqref{eq1ProSta1000star2} for each random point $\xi_i$ to obtain the output set $\{y(t,\xi_i)\}_{i=1}^N$.
\item Aggregate results to estimate the statistics of output $y$ as a function of time.
For instance, if we let $z$ denote the $k$-th component of $y=(y_1,\ldots,y_s)$, then its statistical mean $\mathbf{E}^*[z]:\mathfrak{T}\to\mathbb{R}$ and statistical variance $\mathrm{Var}^*[z]:\mathfrak{T}\to\mathbb{R}^+_0$ are given by:
\begin{equation}
\mathbf{E}^*[z](t)=\frac{1}{N}\sum_{i=1}^N z(t,\xi_i)\quad\text{and}\quad
\mathrm{Var}^*[z](t)=\frac{1}{N-1}\sum_{i=1}^N \big(z(t,\xi_i)-\mathbf{E}^*[z](t)\big)^2.
\end{equation}
\end{enumerate}

\end{appendices}

\section*{Acknowledgements}
The first author gratefully acknowledges the scholarship granted by Colciencias and Atl\'antico Department (Colombia) under Call 673 to pursue a Ph.D.~degree in Structural Engineering at Purdue University, West Lafayette, IN.
The authors also gratefully acknowledge the support of the National Science Foundation (CNS-1136075, DMS-1555072, DMS-1736364, CMMI-1634832, and CMMI-1560834), Brookhaven National Laboratory subcontract 382247, ARO/MURI grant W911NF-15-1-0562, and U.S.~Department of Energy (DOE) Office of Science Advanced Scientific Computing Research program DE-SC0021142.

\section*{Credit author statement}
{\bf H.E.}~conceived the mathematical models, implemented the methods, designed the numerical experiments, interpreted the results, and wrote the paper. {\bf A.P.}~and {\bf G.L.}~supported the study, edited, and reviewed the final manuscript. All the authors gave their final approval for publication.

\bibliography{references1}

\end{document}